\documentclass[11pt]{amsart}

\usepackage[english]{babel}
\usepackage{amssymb,amsfonts}
\usepackage{mathtools}
\usepackage{graphicx}
\usepackage[letterpaper,top=2cm,bottom=2cm,left=3cm,right=3cm,marginparwidth=1.75cm]{geometry}
\usepackage[colorlinks=true,allcolors=blue]{hyperref}
\usepackage{multirow}
\usepackage[nameinlink,noabbrev]{cleveref}
\usepackage{float}

\title[Multiple Timescale Dynamics of Locomotor Neurons]{Multiple Timescale Dynamics of Conductance-based Models of Brainstem Locomotor Neurons}
\date{} 

\author[A. K. Thomas]{Anna Kishida Thomas$^{\dagger}$}
\author[J. E. Rubin]{Jonathan E. Rubin$^{\dagger}$}

\thanks{$^{\dagger}$Department of Mathematics, University of Pittsburgh; and Center for the Neural Basis of Cognition, Pittsburgh, PA, USA.
\texttt{akt51@pitt.edu}, \texttt{jonrubin@pitt.edu}.}

\begin{document}

\begin{abstract}
The pedunculopontine nucleus (PPN) is a heterogeneous brainstem locomotor hub implicated in Parkinson’s disease and potentially relevant for its treatment. We propose single-compartment, conductance-based models for three classes of PPN neurons, such that each model reproduces relevant experimentally observed stimulus-dependent responses, including post-inhibitory rebound dynamics, transient low-threshold activity, and gamma band oscillations. To understand the mechanisms underlying these transient responses to current stimulation, we leverage the models' intrinsic multi-timescale structure and apply dynamical system methods designed for multiple timescale systems. By separating fast membrane and channel-gating dynamics from slower gating and calcium processes, we identify specific ionic mechanisms underlying hallmark dynamics across cell types. We also generate new predictions about PPN behavior under a post-inhibitory facilitation protocol.
\vspace{0.2cm}

\noindent \textit{Relevance to Life Science}: This work links distinct PPN firing responses to specific ion-channel mechanisms.
We develop conductance-based models for three classes of PPN neurons: cholinergic (C), cholinergic with low-threshold calcium spikes (CT), and noncholinergic (NC). On a mechanistic level, we identify post-inhibitory delays in the C model to arise from slow A-current gating; in the CT model, post-inhibitory rebound oscillations are controlled primarily by T-type calcium kinetics; and in the NC model, P/Q- and T-type calcium conductances shape high-voltage oscillations and low-threshold spikes. The models further predict distinct forms of post-inhibitory facilitation selectively in CT and NC classes that express low-threshold calcium currents. Together, these results suggest how inhibitory signals, such as from the largely GABAergic substantia nigra pars reticulata, and DBS-like inputs may differentially shape output across PPN neuron types in brainstem motor circuits.
\vspace{0.2cm}

\noindent \textit{Mathematical Content}: We apply a multiple-timescale dynamical systems framework to analyze neuronal models representing various PPN neuron types. In this approach, we nondimensionalize the models to decompose model variables into fast, slow, and (in the CT model) superslow classes. Within this framework, we apply geometric singular perturbation theory (GSPT) to analyze model dynamics in terms of interactions of reduced subsystems that lead to a variety of bifurcations. In particular, we explain experimentally observed, stimulus-induced phenomena such as rebound oscillations, delays, and transient spiking. Distinct from most previous analyses of this type, we find that capturing some of these responses in the GSPT framework requires taking into account the voltage-dependence of the timescales of certain gating variables.
\end{abstract}

\maketitle




\section{Introduction}

\noindent The pedunculopontine nucleus (PPN) is a brainstem structure within the mesencephalic locomotor region (MLR) where many signaling pathways converge. The PPN plays a central role in numerous brain functions including sleep-wake regulation, initiation and cessation of movement, respiration, and valence processing~\cite{hyam2019,lima2019,fallah2024inhibitory,lin2023ppn,kroeger2017cholinergic}. It is increasingly recognized as relevant for Parkinson's disease (PD) and serves as a target for deep brain stimulation (DBS) in treatment of PD motor deficits, as it is extensively connected to other subcortical brain regions along various motor pathways~\cite{lin2023ppn}. Strategically positioned near the output nuclei of the basal ganglia (BG), the PPN forms direct synaptic connections with the substantia nigra pars reticulata (SNr), subthalamic nucleus (STN), thalamus, medulla and spinal cord~\cite{lee2000pedunculopontine}. Recent studies have demonstrated that transient stimulation of specific BG sites can promote long-term motor recovery in dopamine-depleted rodents, but notably only when the PPN remains intact~\cite{cundiff2024}.

Experiments have revealed different classes of PPN neurons, highlighting the heterogeneity of this nucleus and its diverse dynamics~\cite{kang1990electrophysiological}. These classes include cholinergic (C), cholinergic with T-type calcium current (CT), and noncholinergic (NC) PPN cells, distinguished by their major neurotransmitter and the existence or not of characteristic low-threshold calcium spikes (LTS)~\cite{fallah2024inhibitory,wang2009pedunculopontine, petzold2015decoding, kang1990electrophysiological,takakusaki1997ionic,luster2016intracellular,kezunovic2011mechanism}. Neurons exhibiting the LTS property can present transient spiking activity after the application of an inhibitory current, attributed to low-voltage activated calcium channels. To investigate the factors hypothesized to contribute to certain PPN firing patterns under control and PD conditions, we can turn to computational models of PPN neurons. In particular, conductance-based models explicitly represent transmembrane ion channels and their gating kinetics, enabling direct testing of mechanisms underlying observed behaviors.
In previous work, analyses leveraging this framework have uncovered key mechanisms involved in forms of neuronal voltage dynamics such as burst firing and mixed-mode oscillations~\cite{coombes2005,desroches2012,bertram2017}.
Generic neuronal models of this type cannot capture observations specific to PPN neurons, however.
Moreover, existing models of PPN neurons either consist of general firing rate equations developed to study the population-level activity of the PPN and its impact on other brain regions~\cite{hadipour2004effects,zhao2024dynamic,hu2026multimodal}, focus largely on PPN morphology \cite{zitella2013computational}, or include small sets of currents initially tuned using experiments done outside of PPN and calibrated to a limited range of dynamics in a single PPN cell type \cite{lourens2011pedunculopontine}. These models omit characteristic high-voltage currents and fail to reproduce many of the distinctive forms of single-cell PPN dynamics observed experimentally. 

We address this gap by developing three Hodgkin--Huxley--type conductance-based models from the ground up. Each current in the three models plays an important role in shaping the PPN neuronal activity that has been observed in response to applied current patterns. The models introduced in this work reliably capture diverse post-stimulus rebound responses, gamma band (30$-$90~Hz) oscillations, and transient low-threshold activity seen in their respective cell types. In addition, we predict PPN responses to an untested current protocol for post-inhibitory facilitation \cite{dodla2006,rubin2021}. To systematically investigate and analyze the transient model activity during stimulus-induced responses, we exploit the presence of multiple timescales within the model systems. The slow timescales of gating mechanisms for low-threshold calcium and potassium channels and intracellular calcium dynamics, relative to the fast changes in the transmembrane potential and other gating variables, allow for the application of geometric singular perturbation theory (GSPT). Within this framework, a simplification of the model system produces a fast subsystem, which is studied to understand the transient oscillations and delays arising in the overarching system. GSPT also allows for the isolated study of the slow subsystem to investigate its longer-term contributions to firing patterns. The dependence of the fast subsystem on the slower variables and the dependence of full single-cell dynamics on the pattern of applied current can be studied using bifurcation analysis. Our approach emphasizes the intrinsic dynamics of PPN neurons, laying a foundation for future investigations into PPN function within PD-related motor circuits.

Overall, the objective of this paper is to explain single-cell activity of PPN neurons by identifying key ionic current mechanisms in models for three PPN neuron types. Using the NC model, we show how the conductance of both P/Q- and T-type calcium currents shape the high-voltage oscillations and low-threshold calcium spikes seen in this cell type. After a transient inhibitory current is applied, the C model exhibits a delay before returning to regular firing, which we show relates to the slow gating kinetics of the A-current. The CT model matches experimentally observed rebound oscillations, where multiple timescale analysis highlights the T-type calcium kinetics as the most important factor in shaping the response. Lastly, we test a new post-inhibitory facilitation (PIF) protocol in our models and derive the prediction that only the two cell types containing T-type calcium currents (CT and NC) exhibit PIF, each with unique features. The remaining sections of this manuscript are organized as follows. The proposed PPN models are introduced in~\Cref{sec:models}, along with an overview of the nondimensionalization procedure used in multiple timescale analysis. Beginning with a core model containing a few fundamental currents that produce spiking dynamics, each subsection of the results in~\Cref{sec:results} introduces the contributions to the full models of a few ionic currents at a time, and each result lays out the mechanistic underpinnings driving a specific voltage response. The reader is guided through a computational analysis of dynamics with biological interpretations of the contributing factors, followed by a mathematical analysis explaining the intrinsic dynamics mechanistically.
Conclusions and discussion follow in~\Cref{sec:discussion}, and the appendix offers a brief introduction to our approach to multiple timescale analysis for non-specialist readers and a detailed description of the nondimensionalization procedure illustrated using the CT model.

\section{Conductance-based models for PPN voltage dynamics}
\label{sec:models}

\subsection{Model formulation and parameter values}
We model three types of PPN neurons, distinguished by their major neurotransmitter and characteristic ionic currents: C, CT, and NC. The equations for each of these model types, which describe a single-compartment neuron based on the Hodgkin-Huxley framework \cite{hodgkin1952quantitative}, are
\begin{subequations}
\begin{align}
    \frac{dV_{j}}{dt} &= -\sum_{i \in \mathcal{I}_{j}} I_i + I_{\mathrm{App}} \quad \mathrm{where} \; j \in \{\mathrm{C,CT,NC}\},\label{SYS:PPN:DVDT}\\
    \frac{dp}{dt} &= \frac{p_{\infty}(V_{j}) - p}{t_p(V_{j})},\label{SYS:PPN:DPDT}\\
    \frac{d\mathrm{Ca}}{dt} &= f_{\mathrm{Ca}} \left(-\frac{\sum_{i \in \mathcal{I}_{j,\mathrm{Ca}}} I_i}{2 \cdot F \cdot \mathrm{\textit{Vol}}} - \frac{\mathrm{Ca}-\mathrm{Ca}_{\textit{eq}}}{t_{\textit{store}}}\right), \label{SYS:PPN:DCDT}
\end{align}
\label{SYS:PPN}
\end{subequations}

\noindent where $V_j$ is the neuronal membrane potential for cell type $j \in \{\mathrm{C,CT,NC}\}$, $t$ is the time variable, Ca is the intracellular Ca$^{2+}$ concentration, and each $p \in \{m_{\mathrm{Na}}$, $h_{\mathrm{Na}}$, $m_\mathrm{K}$, $m_{\mathrm{CaPQ}}$, $m_{\mathrm{A}}$, $h_{\mathrm{A}}$, $m_{\mathrm{CaT}}$, $h_{\mathrm{CaT}}\}$ represents a voltage-dependent gating variable, with activation and deactivation gates $m_i$ and $h_i$, respectively, for current $I_i$. We omit capacitance from Eq.~\eqref{SYS:PPN:DVDT} and specify maximal conductances per unit capacitance, rather than assuming capacitance values for the models. For each cell type, we define the set of calcium currents $\mathcal{I}_{j,\mathrm{Ca}}$ contributing to intracellular calcium concentration in Eqs.~\eqref{EQN:CHOL},~\eqref{EQN:CHOLT},~\eqref{EQN:NONCHOL} below.
Note that the subscript of $V_j$ is not strictly necessary, since we instantiate the full set of Eq.~\eqref{SYS:PPN} separately for each neuron type, but we introduce it here so that we can use it to label simulation results later in the paper.

The right-hand side of the first Eq.~\eqref{SYS:PPN:DVDT} consists of the sum of ionic currents governing the dynamics of the membrane voltage for whichever cell type is being modeled. In addition to the standard fast Na$^+$ current $(I_{\mathrm{Na}})$, fast K$^+$ current $(I_{\mathrm{K}})$, and leakage current $(I_{\mathrm{L}})$, the models include contributions from a variety of other currents identified experimentally in PPN cells, which we denote as follows:
\begin{subequations}
\begin{align}
\mathcal{I}_{\mathrm{C}} &= \{I_{\mathrm{Na}}, I_{\mathrm{K}}, I_{\mathrm{L}}, I_{\mathrm{CaPQ}}, I_{\mathrm{KCa}}, I_{\mathrm{A}}\},
&\; \mathcal{I}_{\mathrm{C,Ca}} &= \{I_{\mathrm{CaPQ}}\},
\label{EQN:CHOL}\\
\mathcal{I}_{\mathrm{CT}} &= \{I_{\mathrm{Na}}, I_{\mathrm{K}}, I_{\mathrm{L}}, I_{\mathrm{CaPQ}}, I_{\mathrm{KCa}}, I_{\mathrm{A}}, I_{\mathrm{CaT}}\},
&\; \mathcal{I}_{\mathrm{CT,Ca}} &= \{I_{\mathrm{CaPQ}}, I_{\mathrm{CaT}}\},
\label{EQN:CHOLT}\\
\mathcal{I}_{\mathrm{NC}} &= \{I_{\mathrm{Na}}, I_{\mathrm{K}}, I_{\mathrm{L}}, I_{\mathrm{CaPQ}}, I_{\mathrm{KCa}}, I_{\mathrm{CaT}}\},
&\; \mathcal{I}_{\mathrm{NC,Ca}} &= \{I_{\mathrm{CaPQ}}, I_{\mathrm{CaT}}\},
\label{EQN:NONCHOL}
\end{align}
\end{subequations}

\noindent for C, CT, and NC neurons, respectively. Specifically, the additional currents are high-voltage activated P/Q-type Ca$^{2+}$ current $(I_{\mathrm{CaPQ}})$, Ca$^{2+}$-activated K$^+$ current $(I_{\mathrm{KCa}})$, A-type K$^+$ current $(I_{\mathrm{A}})$, and low-voltage activated T-type Ca$^{2+}$ current $(I_{\mathrm{CaT}})$. $I_{\mathrm{App}}$ denotes an applied current injected from an electrode. The equations for currents follow the Hodgkin-Huxley formulation~\cite{hodgkin1952quantitative}
\begin{equation}
\label{sys:curr}
\begin{aligned}
    I_{\mathrm{Na}} (V_{j}, m_{\mathrm{Na}}, h_{\mathrm{Na}}) &= g_{\mathrm{Na}} \cdot m_{\mathrm{Na}}^3 \cdot h_{\mathrm{Na}} \cdot (V_{j} - E_{\mathrm{Na}}),\\
    I_{\mathrm{K}} (V_{j}, m_{\mathrm{K}}) &= g_{\mathrm{K}} \cdot m_{\mathrm{K}}^4 \cdot (V_{j} - E_{\mathrm{K}}),\\
    I_{\mathrm{L}}(V_{j}) &= g_{\mathrm{L}} \cdot (V_{j} - E_{\mathrm{L}}),\\
    I_{\mathrm{CaPQ}} (V_{j}, m_{\mathrm{CaPQ}}) &= g_{\mathrm{CaPQ}} \cdot m_{\mathrm{CaPQ}} \cdot (V_{j} - E_{\mathrm{Ca}}),\\
    I_{\mathrm{KCa}} (V_{j}, \mathrm{Ca}) &= g_{\mathrm{KCa}} \cdot m_{\mathrm{KCa}}(\mathrm{Ca}) \cdot (V_{j} - E_{\mathrm{K}}),\\
    I_{\mathrm{A}} (V_{j}, m_{\mathrm{A}}, h_{\mathrm{A}}) &= g_{\mathrm{A}} \cdot m_{\mathrm{A}} \cdot h_{\mathrm{A}} \cdot (V_{j} - E_{\mathrm{K}}),\\
    I_{\mathrm{CaT}} (V_{j}, m_{\mathrm{CaT}}, h_{\mathrm{CaT}}) &= g_{\mathrm{CaT}} \cdot m_{\mathrm{CaT}}^2 \cdot h_{\mathrm{CaT}} \cdot (V_{j} - E_{\mathrm{Ca}}),
\end{aligned}
\end{equation}
\noindent where $g_i$ is the maximum conductance and $E_i$ is the reversal potential for each current $I_i$. 

The second Eq.~\eqref{SYS:PPN:DPDT} describes the dynamics governing each activation or inactivation gating variable $p$. The voltage-dependent steady-state function and voltage-dependent relaxation time of $p$ are described by
\[
p_{\infty}(V_{j}) = \frac{1}{1 + \exp\!\left(-(V_{j}-p_{1/2})/k_p\right)},
\]
\[
t_p(V_{j}) = t_p^0 + \frac{t_p^1-t_p^0}{\exp\!\left((\theta_p-V_{j})/\sigma_p^0\right) + \exp\!\left((\theta_p-V_{j})/\sigma_p^1\right)}.
\]

\noindent The half-activation value of the steady-state function occurs at $p_{1/2}$ with slope $k_p$. An exception to these gating dynamics arises in the case of the calcium-gated potassium channels, where the activation $m_{\mathrm{KCa}}$ follows the dynamics outlined in~\cite{benison2001modeling},
\[
    m_{\mathrm{KCa}}(\mathrm{Ca}) = \frac{\mathrm{Ca}^4}{\mathrm{Ca}^4 + K_{\mathrm{KCa}}^4},
\]
\noindent where $K_{\mathrm{KCa}}$ is the affinity constant in the Hill-type activation function, representing the calcium concentration at which the KCa channel is half-activated.

The third Eq.~\eqref{SYS:PPN:DCDT} models the average intracellular calcium concentration, incorporating both the free cytosolic calcium fraction, $f_{\mathrm{Ca}}=0.025$, and a first-order representation of calcium storage dynamics with a time constant $t_{store}=12.5$ ms and baseline concentration $\mathrm{Ca}_{eq}=100$ nM. Intracellular calcium increases due to influx through voltage-dependent channels $\mathcal{I}_{\mathrm{j},\mathrm{Ca}}$ in model type $j$. This influx is converted into moles of calcium ions using the Faraday constant and the divalent nature of calcium $(2 \cdot F)$, and the resulting concentration is averaged over the volume of the neuronal compartment $(\textit{Vol}=7.238\times 10^{-6}$ pL) ~\cite{adams2016computational}. PPN neurons have been shown to exhibit spike-frequency adaptation (SFA) in single-cell electrophysiology experiments \cite{simon2010gamma}. SFA usually results from slowly activating or deactivating currents such as calcium-dependent potassium currents. We implemented $I_{\mathrm{KCa}}$ in all three PPN model types to match the observations of SFA in PPN cells and their concave frequency--applied current $(f-I)$ curves~\cite{simon2010gamma}.

Initial choices of other model parameter values were informed by existing neuronal models (see~\Cref{TAB:PARAMS}) and values were iteratively refined to match the experimental benchmarks described in this section. The final set of parameter values for the models is summarized in~\Cref{TAB:PARAMS}; these represent exemplars of a region in parameter space with dynamically similar dynamics.

\begin{table}[ht!]
\centering
\caption{\label{TAB:PARAMS} Model parameter values for all currents. Conductance values $g_i$ for current $I_i$ in C, CT and NC PPN models, respectively, are separated by commas. Units for current $I_i$ with gating variable $p\in\{m_i,h_i\}$: $g_i~(\mathrm{nS/pF})$, $E_i~(\mathrm{mV})$, $p_{1/2}~(\mathrm{mV})$, $k_p~(\mathrm{mV})$, $t_p^0~(\mathrm{ms})$, $t_p^1~(\mathrm{ms})$, $\theta_p~(\mathrm{mV})$, $\sigma_p^0~(\mathrm{mV})$, $\sigma_p^1~(\mathrm{mV})$, $K_{\mathrm{KCa}}~(\mathrm{nM})$.}
\renewcommand{\arraystretch}{1.3}

\begin{tabular}{lllll}
\hline
Channel & Parameters & & & Reference \\
\hline
$I_{\mathrm{Na}}$ & $g_{\mathrm{Na}} = 50,50,20$ & $E_{\mathrm{Na}} = 50$ & $m_{\mathrm{Na},1/2} = -32.5$ & \cite{phillips2020effects,dautan2021modulation} \\
 & $k_{m_{\mathrm{Na}}} = 7$ & $t^0_{m_{\mathrm{Na}}} = 0.05$ & $t^1_{m_{\mathrm{Na}}} = 0.2$ & \\
 & $\theta_{m_{\mathrm{Na}}} = -12$ & $\sigma^0_{m_{\mathrm{Na}}} = 4$ & $\sigma^1_{m_{\mathrm{Na}}} = -10$ & \\
 & $h_{\mathrm{Na},1/2} = -63$ & $k_{h_{\mathrm{Na}}} = -8$ & $t^0_{h_{\mathrm{Na}}} = 0.7$ & \\
 & $t^1_{h_{\mathrm{Na}}} = 31$ & $\theta_{h_{\mathrm{Na}}} = -48$ & $\sigma^0_{h_{\mathrm{Na}}} = 12$ & \\
 & $\sigma^1_{h_{\mathrm{Na}}} = -6$ & & & \\
\hline
$I_{\mathrm{K}}$ & $g_{\mathrm{K}} = 40,40,0.9$ & $E_{\mathrm{K}} = -87$ & $m_{\mathrm{K},1/2} = -32$ & \cite{traub1992computer,phillips2020effects} \\
 & $k_{m_{\mathrm{K}}} = 10$ & $t^0_{m_{\mathrm{K}}} = 0.3$ & $t^1_{m_{\mathrm{K}}} = 13$ & \\
 & $\theta_{m_{\mathrm{K}}} = -38$ & $\sigma^0_{m_{\mathrm{K}}} = 19$ & $\sigma^1_{m_{\mathrm{K}}} = -16$ & \\
\hline
$I_{\mathrm{L}}$ & $g_{\mathrm{L}} = 0.1,0.1,0.07$ & $E_{\mathrm{L}} = -59$ & & \cite{phillips2020effects} \\
\hline
$I_{\mathrm{CaPQ}}$ & $g_{\mathrm{CaPQ}} = 0.35,0.35,0.21$ & $m_{\mathrm{CaPQ},1/2} = -23$ & $k_{m_{\mathrm{CaPQ}}} = 5.5$ & \cite{mandge2018biophysically} \\
 & $t^0_{m_{\mathrm{CaPQ}}} = 1$ & $t^1_{m_{\mathrm{CaPQ}}} = 1$ & $\theta_{m_{\mathrm{CaPQ}}} = 1$ & \\
 & $\sigma^0_{m_{\mathrm{CaPQ}}} = 1$ & $\sigma^1_{m_{\mathrm{CaPQ}}} = 1$ & & \\
\hline
$I_{\mathrm{KCa}}$ & $g_{\mathrm{KCa}} = 0.5,0.6,0.5$ & $K_{\mathrm{KCa}} = 400$ & & \cite{benison2001modeling} \\
\hline
$I_{\mathrm{CaT}}$ & $g_{\mathrm{CaT}} = 0,4,1$ & $E_{\mathrm{Ca}} = 60$ & $m_{\mathrm{CaT},1/2} = -53.2$ & \cite{golomb1996propagation} \\
 & $k_{m_{\mathrm{CaT}}} = 6.4$ & $t^0_{m_{\mathrm{CaT}}} = 10$ & $t^1_{m_{\mathrm{CaT}}} = 10$ & \\
 & $\theta_{m_{\mathrm{CaT}}} = 1$ & $\sigma^0_{m_{\mathrm{CaT}}} = 1$ & $\sigma^1_{m_{\mathrm{CaT}}} = 1$ & \\
 & $h_{\mathrm{CaT},1/2} = -76.8$ & $k_{h_{\mathrm{CaT}}} = -4.5$ & $t^0_{h_{\mathrm{CaT}}} = 125$ & \\
 & $t^1_{h_{\mathrm{CaT}}} = 100$ & $\theta_{h_{\mathrm{CaT}}} = 1$ & $\sigma^0_{h_{\mathrm{CaT}}} = 1$ & \\
 & $\sigma^1_{h_{\mathrm{CaT}}} = 1$ & & & \\
\hline
$I_{\mathrm{A}}$ & $g_{\mathrm{A}} = 4,9,0$ & $m_{\mathrm{A},1/2} = -16.5$ & $k_{m_{\mathrm{A}}} = 5.14$ & \cite{li2013two} \\
 & $t^0_{m_{\mathrm{A}}} = 0$ & $t^1_{m_{\mathrm{A}}} = 7.58$ & $\theta_{m_{\mathrm{A}}} = -79$ & \\
 & $\sigma^0_{m_{\mathrm{A}}} = 13.3$ & $\sigma^1_{m_{\mathrm{A}}} = -40.3$ & $h_{\mathrm{A},1/2} = -75.7$ & \\
 & $k_{h_{\mathrm{A}}} = -6$ & $t^0_{h_{\mathrm{A}}} = 0$ & $t^1_{h_{\mathrm{A}}} = 16.82$ & \\
 & $\theta_{h_{\mathrm{A}}} = -104$ & $\sigma^0_{h_{\mathrm{A}}} = 5.1$ & $\sigma^1_{h_{\mathrm{A}}} = -255$ & \\
\hline
\end{tabular}
\end{table}

\subsection{Overview of nondimensionalization procedure}\label{sec:sub:overview-nondim}

The proposed neuron models each consist of eight or more dynamical variables, shaping rich single-cell dynamics. The high dimensionality of the models, however, presents a challenge for analysis. The timescales of the variables within these models differ greatly, and fast-slow decomposition is frequently used to analyze multiple timescale dynamical systems~\cite{bertram2017,park2018multiple,john2024novel}. This framework extracts select, lower-dimensional subsystems that can be studied using one- and two-parameter bifurcation analysis. This section provides an overview of nondimensionalization, a systematic approach to identify groupings of variables by their timescales, for the CT model's post-inhibitory rebound response discussed in~\Cref{sec:ct-pir}. The other result for this model type (see~\Cref{sec:ct-pif}) as well as analyses for the C (see~\Cref{sec:c-delay}) and NC model types (see~\Cref{sec:nc-ramp,sec:nc-sdp,sec:nc-pif}) also involve timescale decomposition; however, the identification of fast and slow subsystems in these cases follows an analogous procedure that is omitted for the sake of brevity. Further details, including the timescales identified via the nondimensionalization procedure for the C and NC models, are summarized in~\Cref{sec:appendix-nondim}. 

We first define new dimensionless variables $(v_{\mathrm{CT}}, \mathrm{ca}, \tau)$. Note that the gating variables in the original model are inherently dimensionless as they represent probabilities of ionic gates being open or closed. Scaling factors $k_v$, $k_{\mathrm{ca}}$, and $k_{\tau}$ for transmembrane voltage, intracellular calcium concentration, and time, respectively, are defined such that

$$V_j = k_v \cdot v_{j}, \;\; \mathrm{Ca} = k_{\mathrm{ca}} \cdot \mathrm{ca}, \;\; t = k_{\tau} \cdot \tau.$$

\noindent For the dynamic PIR response, the voltage range is $V_j \in [-92, 41]$, and we choose $k_v = 100$ mV so that $v_{j}=V_{j}/k_v$ remains between $-$1 and 1. Similarly, we determine $k_{\mathrm{ca}} = 100 \mathrm{nM}$ based on $\mathrm{Ca} \in [100,300]$ and $k_{\tau} = 1000 \mathrm{ms}$.

The goal is to write the governing equation for the dynamics of each nondimensionalized variable $x$ in the CT model in the form
$x' = R_x \Tilde{f_x}$
where $R_x$ is a dimensionless constant such that the function $\Tilde{f}_x$ is $\mathcal{O}(1)$ over relevant ranges of its input arguments. The computations for $R_x$ and formulation of $\Tilde{f}_x$ can be found in Appendix~\ref{sec:appendix-nondim}. Within that framework, numerical evaluation reveals the values for the CT model shown in~\Cref{TAB:ct-nondim}.

Based on~\Cref{TAB:ct-nondim}, we identify $x = [v_{\mathrm{CT}}, m_{\mathrm{Na}}, m_{\mathrm{K}},m_{\mathrm{A}},h_{\mathrm{Na}},m_{\mathrm{CaPQ}}]^{\mathsf{T}}$, $y = [m_{\mathrm{CaT}},h_{\mathrm{A}}]^{\mathsf{T}}$, $z=[h_{\mathrm{CaT}},\mathrm{ca}]^{\mathsf{T}}$ as sets of fast, slow, and superslow variables, respectively. We thus represent the CT model for the PIR response as a fast-slow-superslow system of the form
\begin{equation}\label{sys:ct-nondim}
\begin{aligned}
    \frac{dx}{d\tau} &= f(x,y,z), \\
    \frac{dy}{d\tau} &= \varepsilon_1 g(x,y), \\
    \frac{dz}{d\tau} &= \varepsilon_1 \varepsilon_2 h(x,y,z),
\end{aligned}
\end{equation}
\noindent for fast time $\tau$, with $f,g,h$ and $\varepsilon_1$, $\varepsilon_2$ defined in System~\eqref{sys:appendix-ct-nondim2}. We encourage readers unfamiliar with this type of multiple timescale analysis to see Appendix \ref{apn:GSPT}.

The identification of multiple timescales and corresponding fast, slow, and superslow subsystems by applying nondimensionalization to our PPN models guides many of the analyses presented in this work. For clear comparison between the results of our models and experimental recordings, however, the figures are shown using the original model variables.

\begin{table}[H]\centering
\caption{\label{TAB:ct-nondim} Constants $R_x$ for nondimensionalization and classification of variables in the CT model's PIR (\Cref{sec:ct-pir}) and PIF (\Cref{sec:ct-pif}) responses.}
\renewcommand{\arraystretch}{1.3}
$\begin{array}{l|llllllllll} 
x & v_{\mathrm{CT}} & m_{\mathrm{Na}} & m_{\mathrm{K}} & m_{\mathrm{A}} & h_{\mathrm{Na}} & m_{\mathrm{CaPQ}} & m_{\mathrm{CaT}} & h_{\mathrm{A}} & h_{\mathrm{CaT}} & \mathrm{ca} \\
\hline
R_x \; \left[10^{-3}\right] & 50 & 20 & 2.7 & 2.59 & 1.43 & 1 & 0.1 & 0.1 & 0.008 & 0.001 \\
\mathrm{timescale} & \mathrm{fast} & \mathrm{fast} & \mathrm{fast} & \mathrm{fast} & \mathrm{fast} & \mathrm{fast} & \mathrm{slow} & \mathrm{slow} & \mathrm{superslow} & \mathrm{superslow}
\end{array}$
\end{table}

\section{Explaining model dynamics and underlying ion channel gating mechanisms}
\label{sec:results}

We consider models for three PPN neuron types: cholinergic (C), cholinergic-T (CT), and noncholinergic (NC). These neurons and hence their corresponding models share the same core set of currents: the fast sodium current $I_{\mathrm{Na}}$, the delayed rectifier potassium current $I_{\mathrm{K}}$, and the leakage current $I_{\mathrm{L}}$. These currents come together to form a model comparable to the Hodgkin-Huxley model, with a resting potential near $-60$ mV and with type II excitability, with an abrupt onset of oscillatory spiking at a finite, nonzero frequency when depolarized past a threshold \cite{rinzel2013nonlinear}. To showcase how the models specific to distinct neuron types are carefully developed by adding components to this baseline, we will introduce the remaining currents one-by-one, first highlighting what appears to be the dynamical contribution of each and then providing accompanying mathematical analyses to establish the dynamic mechanisms more precisely. PPN model simulations and the computation of bifurcation diagrams were carried out using XPPAUT~\cite{ermentrout2003simulating} and visualizations were created using XPPLORE~\cite{martin2025xpplore}.\\

\begin{figure}[H]
\centering
\includegraphics{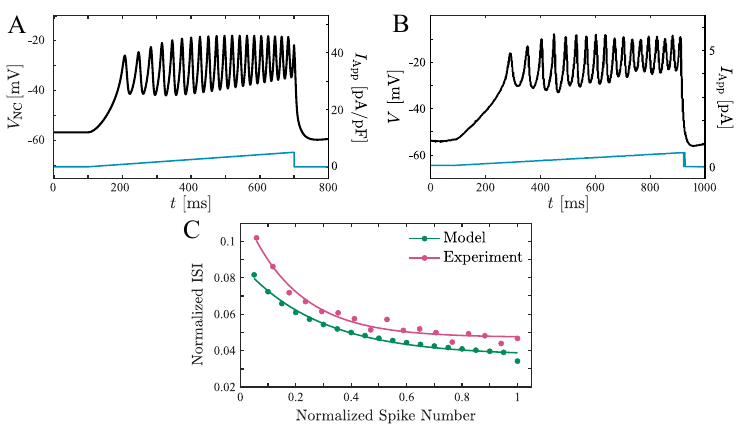}
\caption{\label{fig:nonchol-ramp} \small Comparison of low gamma-band oscillations in response to ramping applied current in the NC model and the experiment by Luster et al.~\cite{luster2016intracellular}.
\textbf{A}: Model voltage trace (black) with ramping current injection $I_{\mathrm{App}}(t), \; t \in [100,700]$ (blue), for $I_{\mathrm{App}}(t)=0.0085\cdot (t-100) \cdot H(t-100) \cdot H(700-t)$ where $H(t)$ is the Heaviside step function. After a gradual rise in voltage, oscillations emerge when $I_{\mathrm{App}}=0.16$ pA/pF.
\textbf{B}: Experimental data reproduced from Fig. 2E (left panel) of Luster et al. ~\cite{luster2016intracellular}, licensed under CC BY 4.0. Voltage traces show intrinsic membrane oscillations during 1-second current ramps from \textit{in vitro} patch-clamp recordings in rat PPN neurons.
\textbf{C}: Comparison of interspike-intervals (ISIs) between NC model (green, see panel A) and experiment~\cite{luster2016intracellular} (pink, see panel B). Oscillations in the model and the experiment both increase in frequency over time as the ramp current is increased. ISIs are normalized to the time window during which the voltage exhibits oscillatory activity. The model and the experimental data (dots) are fit to an exponential function (solid line).}
\end{figure}

\subsection{\texorpdfstring{$I_{\mathrm{CaPQ}}$}{ICaPQ} mediates gamma-band oscillations after transient dynamics involving \texorpdfstring{$I_{\mathrm{CaT}}$}{ICaT} in the NC model during a depolarizing stimulus ramp}\label{sec:nc-ramp}

PPN neurons can sustain tonic firing under strong depolarizing input that would drive many other neuronal types into depolarization block due to insufficient capacity for repolarization. This property has been attributed to the presence of high-voltage-activated (HVA) calcium channels. Luster et al.~\cite{luster2016intracellular} reported that all PPN cells express at least one of two HVA calcium currents, N-type and P/Q-type, and that some P/Q-type cells also exhibit LTS properties. Since Hodgkin-Huxley-style formulations for N- and P/Q-type channels share similar activation kinetics and differ primarily in their inactivation, we focus on the P/Q-type current for simplicity. Accordingly, the P/Q-type calcium current $I_{\mathrm{CaPQ}}$ is incorporated into all three model variants alongside sodium, potassium, and leak currents from the core model. We examine the role of $I_{\mathrm{CaPQ}}$ in supporting low gamma-band oscillations using experimental protocols designed to isolate the HVA current, focusing on the NC model, which exhibits LTS mediated by $I_{\mathrm{CaT}}$ that we shall see contributes to the experimentally observed dynamics.

In experiments by Luster et al.~\cite{luster2016intracellular}, tetrodotoxin (TTX) was applied to block voltage-gated sodium channels and isolate HVA calcium channel contributions in rat PPN neurons. The protocol~(\Cref{fig:nonchol-ramp}B) holds the membrane potential between $-50$ and $-60$ mV, a subthreshold level chosen to prevent spiking activity, before applying a nearly one-second linear current ramp. These steps produce gradual depolarization followed by oscillations. The model reproduces both, including a slight increase in oscillation frequency during the ramp~(\Cref{fig:nonchol-ramp}A). Comparison of normalized inter-spike intervals (\Cref{fig:nonchol-ramp}C) shows close agreement between simulation and experiment, with initially longer intervals that decrease approximately monotonically over time within the range considered to represent gamma oscillations. With sodium currents blocked, large spikes are not observed; instead, persistent oscillations of approximately 20 mV amplitude are sustained across a broad current range. At higher voltages, activation of P/Q-type calcium channels enhances depolarization, which recruits outward potassium currents to repolarize the membrane and maintain rhythmic activity~\cite{luster2016intracellular}. 

\begin{figure}[H]
\centering
\includegraphics{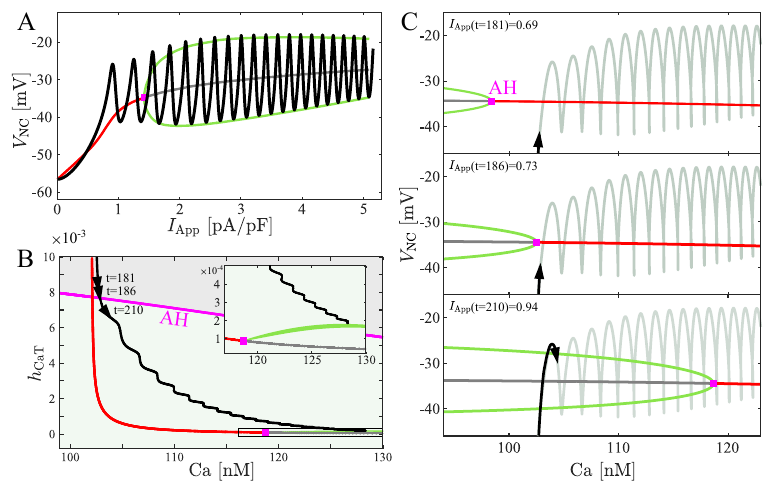}
\caption{\label{fig:nc-ramp-analysis} \small Slow dynamics of $I_{\mathrm{CaT}}$ gating variable and intracellular calcium concentration drive oscillation onset in NC model under depolarizing ramp protocol.
\textbf{A}: Voltage trace (black) overlaid with the bifurcation diagram of the full NC model computed with respect to $I_{\mathrm{App}}$. A branch of stable equilibria (red) persists until $I_{\mathrm{App}}=1.41$, where a supercritical AH bifurcation (pink) gives rise to stable limit cycles (green, minima/maxima). The voltage follows the stable equilibrium branch initially, with oscillations emerging before the AH point. After the bifurcation, $v_{\mathrm{NC}}$ remains oscillatory and generally aligns with the limit cycle extrema.
$\textbf{B}$: Same trajectory (black), full system stable equilibria (red), and family of stable periodic orbits (POs) (green) as in \textit{A}, projected onto the $(\mathrm{Ca}$, $h_{\mathrm{CaT}})$-plane. The trajectory initially tracks the stable equilibrium branch before crossing a branch of fast subsystem AH bifurcations (pink) as applied current increases, moving from a region where the fast subsystem converges to stable equilibria (gray) into a region where the fast system follows stable POs (light green). The inset shows a zoomed view of the boxed region in the main diagram, where the full system AH bifurcation point from \textit{A} occurs.
\textbf{C}: Voltage traces (black) are shown for $t\leq 181$ (top), $\leq 186$ (center), and $\leq 210$ (bottom), overlaid with the $\mathcal{M}_1^{181}|_{\hat{y}(181)}$, $\mathcal{M}_1^{186}|_{\hat{y}(186)}$, and $\mathcal{M}_1^{210}|_{\hat{y}(210)}$, respectively, for $\hat{y}(t)=[h_{\mathrm{CaT}}(t)]$. The gray-green curve shows the future evolution of the trajectory. Stable (red) and unstable equilibrium (gray) branches are shown in addition to an AH bifurcation (pink) and stable limit cycle extrema (green, minima/maxima). Initially, the fast voltage dynamics drive the trajectory toward a stable fixed point. As applied current increases over time, the AH bifurcation shifts to higher $\mathrm{Ca}$ values and the trajectory follows stable limit cycles.
}
\end{figure}

The ramped applied current varies slowly, suggesting that it may be possible to treat it as a bifurcation parameter.
\Cref{fig:nc-ramp-analysis}A shows the full system trajectory superimposed on the bifurcation diagram of the full model with respect to the parameter $I_{\mathrm{App}}$. Stable periodic orbits (green) emerge from a supercritical Andronov-Hopf (AH) bifurcation along a branch of equilibria (red/black). At the onset of the current ramp at $\tau=0.1$, the membrane potential lies at an $I_{\mathrm{App}}$ equilibrium point. As the current increases, the $v_{\mathrm{NC}}$ coordinate of the trajectory departs from the stable equilibrium branch and exhibits oscillations even before the branch destabilizes at the AH. The discrepancy between the voltage dynamics and the projections of the invariant sets of the full system suggests that a disparity between timescales in the model (see~\Cref{sec:sub:overview-nondim}) may contribute to the dynamics. Despite the slowly ramping nature of the applied current, the static bifurcation diagram is not sufficient to fully explain this time-dependent model response.

To determine how timescale separation shapes the model dynamics, we nondimensionalize the NC PPN model under sodium-channel blocker (\Cref{sec:sub:overview-nondim}, Appendix \ref{sec:appendix-nc-nondim}) to obtain a fast-slow system of the form
\begin{equation}\label{sys:nc-ramp-fasttime}
\begin{aligned}
    \frac{dx}{d\tau} &= f(x,y), \\
    \frac{dy}{d\tau} &= \varepsilon_1 g(x,y),
\end{aligned}
\end{equation}
\noindent with fast $x = [v_{\mathrm{NC}},m_{\mathrm{K}},m_{\mathrm{CaPQ}},m_{\mathrm{CaT}}]^{\mathsf{T}}$and slow variables $y=[h_{\mathrm{CaT}}, \mathrm{ca}]^{\mathsf{T}}$ and fast time $\tau$. The procedures to define functions $f,g,$ and parameter $\varepsilon_1$ on the right-hand side of System~\eqref{sys:nc-ramp-fasttime} are outlined in~\Cref{sec:appendix-nondim}. By taking $\varepsilon_1 \rightarrow 0$ in System~\eqref{sys:nc-ramp-fasttime}, we can fix the slow variables to isolate the \textit{fast subsystem}:

\begin{equation}\label{sys:nc-ramp-fastsub}
\begin{aligned}
    \frac{dx}{d\tau} &= f(x,y), \\
    \frac{dy}{d\tau} &= 0.
\end{aligned}
\end{equation}

The \textit{critical manifold} is defined as the set of equilibrium points of this fast subsystem. We will use the notation 
\begin{equation}
\mathcal{M}_1^m=\left\{(x,y): \; f(x,y)=0, \bar{I}_{\mathrm{App}}=\bar{I}_{\mathrm{App}}^m\right\}, \; m \in \{b,i,e\}
\label{EQN:M1m}
\end{equation}
\noindent to denote the critical manifold with nondimensionalized applied current at constant baseline ($b$), excitatory ($e$), or inhibitory ($i$) levels. (See~\Cref{sec:appendix:ct-nondim} for the steps in the nondimensionalization of currents; note that $\bar{I}_{\mathrm{App}}^m = {I}_{\mathrm{App}}^m/(k_v \cdot g_{\mathrm{max}})$). In this protocol, $I_{App}^b=0$ and during the ramp the applied current varies continuously with time $\tau$; hence, we additionally define
\begin{equation}
\mathcal{M}_1^{\tau} = \left\{(x,y): f(x,y)=0,\ \bar{I}_{\mathrm{App}}=\bar{I}_{\mathrm{App}}(\tau)\right\}.
\label{EQN:Iappt}
\end{equation}

\noindent In figures, we will show projections of the high-dimensional critical manifolds. To be able to represent the manifold, the slow (and superslow, for three-timescale systems that we consider later) variables not represented by an axis in the projection will be fixed at values that we will denote by $\hat{y}$. In the NC model, $y\in\mathbb{R}^2$ and we write $y=(y_i,\hat{y})$ where $y_i$ varies along an axis in a visualization and $\hat{y}$ is fixed. We define
\begin{equation}
\mathcal{M}_1^{\tau}\vert_{\hat{y}(\tau)} = \left\{(x,y) \in \mathbb{R}^{4 \times 2}: f(x,(y_i,\hat{y}(\tau)))=0,\bar{I}_{\mathrm{App}}=\bar{I}_{\mathrm{App}}(\tau) \right\}.
\label{EQN:M1myhat}
\end{equation}

\noindent For visualizations, we analogously define the manifold $\mathcal{M}_1^{t}\vert_{\hat{y}(t)}$ for dimensionalized time $t$. For example, in~\Cref{fig:nc-ramp-analysis}C (top), $\hat{y}(t)=[h_{\mathrm{CaT}}]$, and 
$\mathcal{M}_1^{181}\vert_{\hat{y}(181)}$ denotes the manifold $\mathcal{M}_1^{0.181}$ in nondimensionalized space projected onto dimensionalized $(\mathrm{Ca},V_{\mathrm{NC}})$-space with variable $\hat{y}(t)$ evaluated at $t=181$ ms to draw the manifold. We use dimensionalized time $t$ in the visualizations rather than $\tau$.

Note that in~\Cref{fig:nc-ramp-analysis}B, the slower calcium and $h_{\mathrm{CaT}}$ dynamics diverge from the full system equilibria before $t=181$. 
Leveraging the presence of two timescales within the NC model, in~\Cref{fig:nc-ramp-analysis}C (top), we can see that early in the ramping response, the fast subsystem is approaching a stable fixed point at a depolarized voltage on the projection of $\mathcal{M}_1^{0.181}$, above (with respect to Ca) the AH bifurcation point. As the applied current increases, the fast voltage depolarizes further and the AH bifurcation on $\mathcal{M}_1^{0.186}$ lies at a higher level of intracellular calcium~(\Cref{fig:nc-ramp-analysis}C, center). As $\bar{I}_{\mathrm{App}}$ increases more, we observe the AH bifurcation on $\mathcal{M}_1^{0.210}$ at a much higher calcium level than that of the trajectory. The crossing of the AH bifurcation associated with the fast dynamics results in oscillations in the fast variables of the NC model.

Let $A^{\tau} \subset \mathcal{M}_1^{\tau}$ denote the curve of AH bifurcations on the critical manifold for fixed $\tau$ and define $A^{\tau_{\mathrm{dur}}}$ as the set of AH bifurcations across the duration of the experiment, $\tau_{\mathrm{dur}}=[0,700]$. \Cref{fig:nc-ramp-analysis}B illustrates how $A^{\tau_{\mathrm{dur}}}$ divides the $(h_{\mathrm{CaT}},\mathrm{Ca})$-plane into a region where the fast subsystem exhibits stable oscillations (light green) and a region where the fast system converges to stable critical points (gray). In fact, the trajectory initiates oscillatory activity once it crosses $A^{\tau_{\mathrm{dur}}}$, allowing the P/Q-channel driven oscillations to persist. For the remainder of the trajectory, the system slowly approaches the family of stable periodic orbits (POs) of the full system (\Cref{fig:nc-ramp-analysis}B, inset).

Although the initial holding potential of this experimental protocol was thought to minimize the effects of low-threshold currents, our analysis predicts that the inactivation gating variable $h_{\mathrm{CaT}}$ of $I_{\mathrm{CaT}}$ contributes to the timing of oscillation onset through its slowly evolving dynamics. Similarly, while the Ca-gated $I_{\mathrm{KCa}}$ is not highly active in these conditions, its activation dynamics, mediated by slow changes in calcium, also affects the voltage dynamics. Our results do support the prevailing hypothesis that the high-frequency oscillations that emerge following TTX application and the depolarizing ramp are driven by the dynamic balance between the inward P/Q-type calcium current and the outward potassium current because the AH bifurcation of the full system and its associated stable limit cycles are crucial to enable these oscillations. But our results provide the more nuanced view that while the ability of PPN neurons to produce and sustain firing under large applied currents is mediated by high-voltage-activated calcium channels, other currents may control the timing of oscillation onset in response to gradual current increases.

\subsection{Slow inactivation of T-type calcium channel mediates transient spiking in the NC PPN model}\label{sec:nc-sdp}

A slowly depolarizing potential (SDP) followed by a burst of activity has been observed in response to a subthreshold step of applied current above a low holding potential in PPN cells with LTS~\cite{kang1990electrophysiological}. A plausible scenario is that the T-type calcium conductance, recruited at the holding potential, temporarily remains high while the small excitation provided by the stimulus gradually depolarizes the cell and eventually spikes arise. We now continue with our NC model analysis and show how $I_{\mathrm{CaT}}$ serves as a major driver of this transient spiking response in the NC model. We note that the CT model, with dynamics that also depend on T-type calcium channels, is able to exhibit a transient spiking response to this input protocol as well. However, CT neurons produce other benchmark activity patterns presented in later sections of this paper, which require a larger maximal T-type calcium conductance. Due to this feature, the small-amplitude excitation considered here elicits a greater number of voltage spikes than observed experimentally. Thus, our work suggests that the experimentally observed SDP response was measured in NC PPN cells, and we focus on the NC model to investigate the underlying dynamics.

\begin{figure}[H]
\centering
\includegraphics{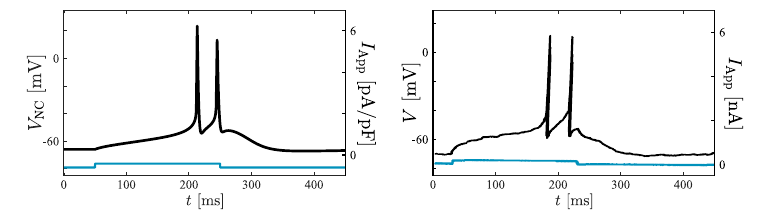}
\caption{\label{fig:nc-sdp} \small Comparison of voltage responses to a small stimulus application in NC model and experimental recordings~\cite{kang1990electrophysiological}. 
\textbf{Left}: Model voltage trace (black) with current step to $I_{\mathrm{App}}(t)=185$ pA/pF (blue). The excitatory stimulus depolarizes the cell, inducing two action potentials, presumed to be mediated by the T-type calcium current, before the injected current is removed and the voltage returns to rest.
\textbf{Right}: Experimental data adapted from Kang and Kitai~\cite{kang1990electrophysiological} in PPN cells exhibiting the LTS property, reproduced with permission from Elsevier. The membrane potential slowly depolarizes after the onset of the stimulus. After approximately $150$ ms, the cell fires two action potentials, after which the stimulus is removed and voltage returns to rest around $-70$ mV.}
\end{figure}

In experiments by Kang and Kitai ~\cite{kang1990electrophysiological}, the membrane potential is held near $-70$ mV before applying a 200 ms subthreshold excitatory pulse (\Cref{fig:nc-sdp}, right). The voltage depolarizes slowly and produces a series of two spikes, then repolarizes after the stimulus is removed. The NC model captures this response (\Cref{fig:nc-sdp}, left). When a small excitatory input is applied from a hyperpolarized state near $-65$ mV, deinactivation of $I_{\mathrm{CaT}}$ permits T-type current recruitment. In response, the voltage slowly depolarizes until it crosses its firing threshold, producing two spikes approximately $30$ ms apart, before the stimulus is removed and the voltage returns back to baseline. The sudden elevation of applied current allows the $m_{\mathrm{CaT}}$ gate to open, while the slow kinetics of the $h_{\mathrm{CaT}}$ inactivation gate prevents $I_{\mathrm{CaT}}$ from deinactivating quickly and hence provides a window of increased $I_{\mathrm{CaT}}$ conductance that represents a natural candidate for the source of the spikes.

The transient subthreshold activity exhibited by the NC model under the SDP protocol can be analyzed in a relatively straightforward way through the lens of a fast-slow timescale decomposition. Through nondimensionalization and multiple timescale identification as outlined in~\Cref{sec:sub:overview-nondim,sec:appendix-nc-nondim}, we arrive at a decomposition into two timescale classes (see~\Cref{tab:nc-timeconstants}) in the NC model during the SDP current protocol: the fast variables $x=[v_{\mathrm{NC}},m_{\mathrm{Na}},m_{\mathrm{K}},\allowbreak h_{\mathrm{Na}},m_{\mathrm{CaPQ}}]^{\mathsf{T}}$ and the slow variables $y=[m_{\mathrm{CaT}},h_{\mathrm{CaT}},\mathrm{ca}]^{\mathsf{T}}$. In the previous subsection, we classified $m_{\mathrm{CaT}}$ as a fast variable; with the introduction of the gating variables for $I_{\mathrm{Na}}$ and the faster evolution of $m_{\mathrm{K}}$ here, $m_{\mathrm{CaT}}$ now acts as a slow variable. The NC model can be represented as a system of the form~\eqref{sys:nc-ramp-fasttime}. The critical manifolds at baseline, $\mathcal{M}_1^b$, and during excitation, $\mathcal{M}_1^e$, are defined using Eq.~\eqref{EQN:M1m} with nondimensionalized $\bar{I}_{\mathrm{App}}^b=-0.6/(k_v \cdot g_{\mathrm{max}})$ and $\bar{I}_{\mathrm{App}}^e=-0.415/(k_v \cdot g_{\mathrm{max}})$, respectively. We will show that the spiking activity in the fast variables is mediated by the slow variables $h_{\mathrm{CaT}}$ and $\mathrm{ca}$.

\begin{figure}[H]
\centering
\includegraphics{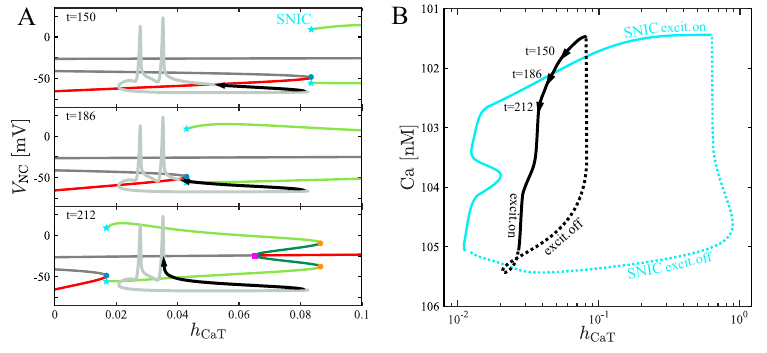}
\caption{\label{fig:nc-sdp-analysis} \small NC model spikes at a fast subsystem saddle-node on an invariant circle (SNIC) bifurcation.
\textbf{A}: Solution to full NC system (black) from time 0 until the time indicated in each panel, along with the future trajectory path (gray-green). Plots include invariant sets of the fast subsystem at $t=150,186,$ and $212$ ms: projections of critical manifolds $\mathcal{M}_1^e\vert_{\hat{y}(150)}$, $\mathcal{M}_1^e\vert_{\hat{y}(186)}$ and $\mathcal{M}_1^e\vert_{\hat{y}(212)}$ for $\hat{y}(t)=[m_{\mathrm{CaT}}(t),\mathrm{ca}(t)]$ (stable red branches and unstable gray branches) and max/min voltages along periodic orbits (green). In each panel, the position of the trajectory at the time indicated is marked with a black triangle. This position aligns with a fast subsystem SNIC bifurcation (cyan) near $t=186$ ms (center panel). After this time, the fast subsystem no longer has a stable branch of the critical manifold nearby, and stable POs become the only attractors of the fast variables. Fast variables follow the PO family for two large-amplitude oscillations as slow variable $h_{\mathrm{CaT}}$ continues to slowly decrease.
\textbf{B}: Illustration of system solution during excitation (black) and after its removals (dotted, black). Black triangles mark the trajectory positions at the indicated times, matching the left panels. Once the trajectory crosses the curve of SNIC bifurcations of $\mathcal{M}_1^e$ with excitation on (solid, cyan), the fast subsystem exhibits two large-amplitude oscillations. After the stimulus is removed, the trajectory ends up above the SNIC curve for $\mathcal{M}_1^b$ (dotted, cyan) and follows a stable branch of equilibria (not explicitly shown) on $\mathcal{M}_1^b$.}
\end{figure}

Before the small-amplitude excitatory step current is applied, the NC model trajectory equilibrates near a full system stable critical point on $\mathcal{M}_1^b$. Once the depolarizing current, $\bar{I}_{\mathrm{App}}^e$, is applied, and before spiking activity is initiated, the trajectory quickly jumps to a surface of stable equilibria on $\mathcal{M}_1^e$ (see snapshots in~\Cref{fig:nc-sdp-analysis}A). Once the fast variables have converged to this attractor, the slow variables, including $h_{\mathrm{CaT}}$, take over the dynamics and drift on the slow timescale $\tau_1$. The slow subsystem moves the fast subsystem along $\mathcal{M}_1^e$ until the fast variables cross a branch of fast subsystem saddle-node on an invariant circle (SNIC) bifurcations. Let $p_{\tau}^{\mathrm{SNIC}} \in \mathcal{M}_1^e\vert_{\hat{y}(\tau)}$ denote the SNIC bifurcation point for fixed $\tau$, then $\mathcal{B}^{\mathrm{SNIC}}:=\cup_{\tau\leq0.25} p_{\tau}^{\mathrm{SNIC}}$ forms a branch of SNIC bifurcations, as seen in~\Cref{fig:nc-sdp-analysis}B. After crossing $\mathcal{B}^{\mathrm{SNIC}}$ near $\tau=0.186$ (\Cref{fig:nc-sdp-analysis}A center), the fast variables move toward the associated family of stable periodic orbits of the fast subsystem. The fast variables complete two oscillation cycles along the family of periodic orbits, as seen in \Cref{fig:nc-sdp-analysis}A (bottom), before the excitatory current is discontinued. The slow drift then takes the fast variables, with $v_{\mathrm{NC}}$ hyperpolarized, back across the SNIC bifurcation (\Cref{fig:nc-sdp-analysis}B). The attractor of the fast dynamics becomes the stable branch of equilibria on $\mathcal{M}_1^b$ at low voltage.
Finally, the slow variables evolve on the slow timescale $\tau_1$ along $\mathcal{M}_1^b$, approaching the attractor of the slow subsystem on this manifold, which is the stable critical point of the full system where the experiment started.

\subsection{A-current leads to delayed return to tonic spiking in the C model}\label{sec:c-delay}

In addition to the high-voltage activated calcium currents, certain PPN cell types express low-voltage

\begin{figure}[H]
\centering
\includegraphics{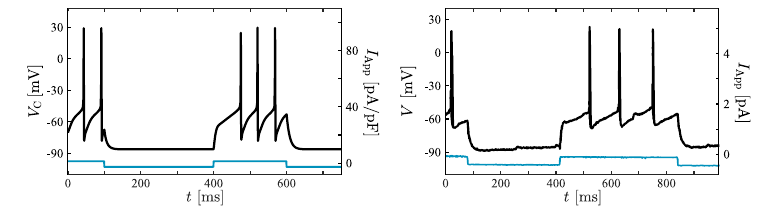}
\caption{\label{fig:chol-delay} \small Comparison of delayed firing after removal of inhibitory input in C model and experimental recordings from rat PPN cells by Takakusaki and Kitai~\cite{takakusaki1997ionic}.
\textbf{Left}: Model voltage trace (black) with inhibitory current (blue) during $t \in (100,400)$ and $(660,800)$. The voltage fires rhythmically before the onset of the inhibition at $t=100$ ms. During the inhibition, the voltage remains hyperpolarized near $-80$ mV. Once the inhibition is removed, the voltage experiences a gradual depolarization before returning to regular firing activity.
\textbf{Right}: Experimental data from Takakusaki and Kitai~\cite{takakusaki1997ionic}, reproduced with permission from Elsevier. The onset of inhibitory current (blue) around $t=90$ ms causes the voltage (black) to remain hyperpolarized near $-90$ mV. After removal of the inhibition near $t=400$ ms, a delayed return to rhythmic firing suggests the presence of the A-current.}
\end{figure}

\noindent activated, transient outward A-type potassium currents.
Indeed, Takakusaki et al.~\cite{takakusaki1997ionic} inferred the presence of a transient outward A-current from intracellular recordings showing a delayed return to regular firing activity upon removal of an inhibitory stimulus in rat PPN cells.
Cells exhibiting $I_{\mathrm{A}}$ are largely immunopositive to choline acetyltransferase; therefore, only our C and CT model variants include this current. To isolate its role, we focus on the C model, which contains $I_{\mathrm{A}}$ without other voltage-dependent low-threshold currents. The C model reproduces the experimentally observed behavior under a comparable protocol (\Cref{fig:chol-delay}). Following the hyperpolarizing step, both experimental and simulated neurons remain below firing threshold for more than 70 ms before the first spike occurs. After this delay, firing resumes at a regular rate, with approximately constant inter-spike intervals (ISIs) and no gradual frequency ramp-up, consistent with the experimentally observed dynamics.

We will illustrate how the A-current is responsible for the characteristic delay in spike initiation following release from hyperpolarization in the C model. Because $I_{\mathrm{A}}$ is activated at low voltages and inactivates slowly, it can transiently oppose depolarization after inhibitory input, shaping the timing of firing recovery. Without this low-threshold current in these cells, the voltage would be expected to rapidly depolarize and return to a tonically firing state. Although this spike-delaying effect of the A-current is well known \cite{gerber1993,ermterm}, we shall analyze the specific bifurcation events through which the slow timescale dynamics of the A-current gating variables contribute to spike onset and timing of successive spikes.

From an analytical perspective, the characteristic delay produced by the A-current can be understood by exploiting the slow timescale of the $I_{\mathrm{A}}$ gating variables $m_{\mathrm{A}}$ and $h_{\mathrm{A}}$. By nondimensionalizing the C model system using the procedure outlined in~\Cref{sec:sub:overview-nondim} and~\Cref{sec:appendix:ct-nondim}, we calculate the time constants for each variable in this model (see~\Cref{tab:c-timeconstants}). Based on these time constants, we group the C model into two subsystems varying on different timescales. The fast subsystem is defined as $x=[v_{\mathrm{C}},m_{\mathrm{Na}},m_{\mathrm{K}},m_{\mathrm{CaPQ}}]^{\mathsf{T}}$ and the slow subsystem is $y=[h_{\mathrm{Na}},m_{\mathrm{A}},h_{\mathrm{A}},\mathrm{ca}]^{\mathsf{T}}$. The C model does not include $I_{\mathrm{CaT}}$, so its gating variables from our previous NC model analysis are not present here. Although the components of $x$ and $y$ differ between the two decompositions, the C model system can be expressed in the same form as the NC model, given by System~\eqref{sys:nc-ramp-fasttime}. Following the procedure to define critical manifolds of fast variables $x$ in the two-timescale system of the NC model (see Eq.~\eqref{EQN:M1m}), we define $\mathcal{M}_1^b$ and $\mathcal{M}_1^i$ as the critical manifolds for baseline $\bar{I}_{\mathrm{App}}^b$ and inhibitory $\bar{I}_{\mathrm{App}}^i$ applied current values.

\begin{figure}[H]
\centering
\includegraphics{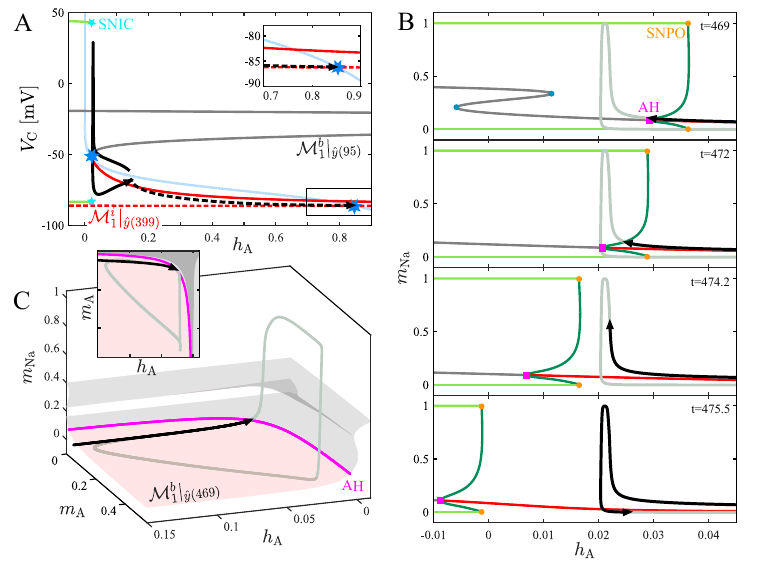}
\caption{\label{fig:chol-delay-analysis} \small Delayed return to firing after removal of an inhibitory stimulus in the C model.
\textbf{A}: The trajectory before (during) the application of inhibition is shown in solid black (dashed), corresponding to voltage from \Cref{fig:chol-delay} (left) for $t$ between $58$ and $400$ ms, when inhibition ends. Black triangles denote points along trajectory at $t=58$ and $400$ ms.
Additional curves are the $h_{\mathrm{A}}$-nullcline (pale blue), the projections $\mathcal{M}_1^b\vert_{\hat{y}(95)},\hat{y}(t)=[h_{\mathrm{Na}}(t),m_{\mathrm{A}}(t),\mathrm{ca}(t)]$ (solid red and gray) and $\mathcal{M}_1^i\vert_{\hat{y}(399)}$ (dashed red), and stable POs of the fast subsystem (max and min of dimensionalized $V_{\mathrm{C}}$ in green).
Before the inhibitory input, the full system has an equilibrium point (blue star near $(h_{\mathrm{A}},V_{\mathrm{C}})=(0.03,-52)$) and is in an oscillatory regime. With the inhibitory input $\bar{I}_{\mathrm{App}}^i$, the model trajectory (black dashed) converges to a stable equilibrium point (blue star near $(h_{\mathrm{A}},V_{\mathrm{C}})=(0.85,-86)$) where the projection $\mathcal{M}_1^i\vert_{\hat{y}(399)}$ and the nullcline of $h_{\mathrm{A}}$ (pale blue) intersect. 
The light blue stars denote a SNIC bifurcation of the fast subsystem. The inset shows the intersection of $\mathcal{M}_1^i\vert_{\hat{y}(399)}$ and $h_{\mathrm{A}}$-nullcline at the blue star.
\textbf{B}: Solution to full C model (black) from time $460$ until the time indicated in each panel, along with the future trajectory path (gray-green) until time $477$ ms. Plots include the invariant sets of the fast subsystem. Projections are shown of critical manifolds $\mathcal{M}_1^b\vert_{\hat{y}(t)}$ for $\hat{y}(t)=[h_{\mathrm{Na}}(t),m_{\mathrm{A}}(t),\mathrm{ca}(t)]$ (stable red branches and unstable gray branches) and max/min voltages along periodic orbits (green) at $t=469, 472, 474.2,$ and $475.5$ ms. In each panel, the position of the trajectory at the time indicated is marked with a black triangle. This position aligns with a fast subsystem AH bifurcation (pink) near $t=469$ ms (top panel). After this time, the fast subsystem jumps to stable periodic orbits of the fast variables. Fast variables follow the PO family and, as slow variables drift, fast variables are attracted to a stable branch of equilibria (bottom panel) during the hyperpolarization phase of each oscillation.
\textbf{C}: Voltage trace from~\Cref{fig:chol-delay} (left) from $t \in (436, 469)$ (black) and $(469,500)$ (gray-green) shown in $(h_{\mathrm{A}}, m_{\mathrm{A}}, m_{\mathrm{Na}})$-space. $\mathcal{M}_1^b\vert_{\hat{y}(469)}$ and the projection of $A^{\tau}$ (pink) into dimensionalized space on this manifold are shown at baseline for fixed $\hat{y}(t)=[h_{\mathrm{Na}}(t),\mathrm{ca}(t)]$ at $t=469$. The subcritical AH of the fast subsystem forms a boundary between stable (pink) and unstable (gray) subregions of $\mathcal{M}_1^b$; around $t=469$, the trajectory passes through this boundary. The inset shows the fast jump away from $\mathcal{M}_1^b$ near the crossing of the AH. The black triangle marks the system state along the trajectory at $t=469$ ms.}
\end{figure}

In the experimental protocol, before any inhibition is applied, the neuronal model is subject to a small input ($\bar{I}_{\mathrm{App}}^b=1.3/(k_v \cdot g_{\mathrm{max}})$, nondimensionalized) and is in a tonic spiking regime (\Cref{fig:chol-delay}, left). \Cref{fig:chol-delay-analysis}A shows the intersection of the $h_{\mathrm{A}}-$nullcline and the projection of $\mathcal{M}_1^b$ at small $h_{\mathrm{A}}$, where there are stable periodic orbits, corresponding to the baseline spiking in the model. From time $\tau=0.1$ to $0.4$, an inhibitory step current ($\bar{I}_{\mathrm{App}}^i=-2.7/(k_v \cdot g_{\mathrm{max}})$) is given and the system settles to a stable equilibrium near $(h_{\mathrm{A}},v_{\mathrm{C}})=(0.85,-0.86)$ where $\mathcal{M}_1^i$ and the $h_{\mathrm{A}}$-nullcline intersect, and the fast subsystem has no periodic orbits. \Cref{fig:chol-delay-analysis}A shows the projection of $\mathcal{M}_1^i\vert_{\hat{y}(\mathrm{399})}$ for $\hat{y}(t)=[h_{\mathrm{Na}},m_{\mathrm{A}},\mathrm{ca}(t)]$. We will show that, once inhibition is removed, there is a delay while we wait for $h_{\mathrm{A}}$ to decrease such that the oscillations can resume.

Once the stimulus is removed, the relevant critical manifold is once again $\mathcal{M}_1^b$, which has a branch of stable equilibria at a higher voltage value. The fast variables approach this branch on the fast timescale $\tau$, after which drift in the slow variables on the slow timescale $\tau_1$ takes over (\Cref{fig:chol-delay-analysis}B, $t$=469).

As the slow variables drift, the trajectory evolves along the lower branch $\mathcal{M}_1^b$ of the S-shaped critical manifold, where a surface of subcritical AH bifurcations emerges (via a Takens--Bogdanov bifurcation, not shown). Let $A^{\tau} \subset \mathcal{M}_1^{\tau}$ denote the family of AH bifurcations on the critical manifold for fixed $\tau$. The system remains near $A^{\tau}$ until, around $469$ ms, the fast variables depart from $\mathcal{M}_1^b$ (see \Cref{fig:chol-delay-analysis}B-C) and approach a stable family of fast subsystem

\noindent POs, leading to a fast oscillation. During this event, the slow variables evolve sufficiently that the trajectory crosses back through the AH surface, yielding a return to a stable part of $\mathcal{M}_1^b$ after the oscillation is complete (\Cref{fig:chol-delay-analysis}B, $t$=472,474.2). Specifically, after the large-amplitude excursion, there is no longer bistability between POs and equilibria of the fast subsystem. Thus, during the afterhyperpolarization phase (\Cref{fig:chol-delay-analysis}B, $t$=475.5), the fast variables settle on $\mathcal{M}_1^b$ and slow drift resumes, with $h_{\mathrm{A}}$ increasing. Once $h_{\mathrm{A}} \approx 0.14$ (\Cref{fig:chol-delay-analysis}C), $h_{\mathrm{A}}$ decreases, initiating another gradual depolarization toward the next AH crossing and spike. 

\vspace{0.1cm}
\textit{\textbf{Remark}}: During the slow depolarization preceding each spike, the trajectory briefly enters the unstable subregion of $\mathcal{M}_1^b$ beyond the AH surface without immediately departing from it. In this region, the positive real parts of the linearization about $\mathcal{M}_1^b$ remain small. Although we do not fully resolve this transient behavior, these weak instabilities likely delay the repulsion from the unstable manifold, allowing the trajectory to remain near $\mathcal{M}_1^b$ before the fast jump occurs.
\vspace{0.1cm}

\subsection{Slow gating of the T-type calcium current induces post-inhibitory stimulus oscillations in the CT model}\label{sec:ct-pir}

Synaptic connectivity between the substantia nigra pars reticulata (SNr), the primary basal ganglia output nucleus in rodents, and the PPN is increasingly receiving attention \cite{fallah2024inhibitory,falasconi2025}, as it links two key brain areas along the motor pathway. SNr neurons are predominantly GABAergic, inducing inhibitory currents that may suppress firing in their synaptic targets. However, some PPN neurons exhibit post-inhibitory rebound (PIR), in which activity increases following the removal of inhibition. In the experiments considered here, this rebound activity is preceded by a gradual depolarization that delays spike onset, distinguishing it from more abrupt PIR phenomena. To explain this behavior, we analyze the CT model, which includes both $I_{\mathrm{A}}$ and a large maximal $I_{\mathrm{CaT}}$, and emphasize the interaction of fast, slow, and superslow gating variables. Unlike the NC and C models, which contain only one low-threshold current, the CT model captures the delayed rebound followed by sustained spiking, highlighting the necessity of both A-type potassium and T-type calcium conductances for this response. 

\begin{figure}[H]
\centering
\includegraphics{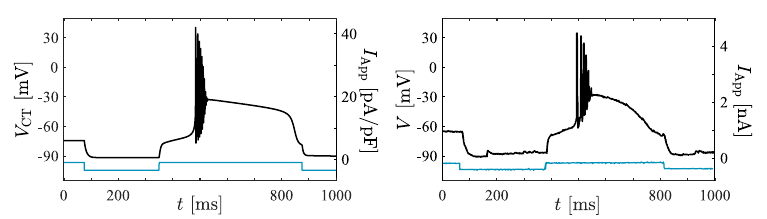}
\caption{\label{fig:cholt-pir} \small Comparison of rebound response to inhibitory step current in the CT model and experiment in Takakusaki et al.~\cite{takakusaki1997ionic}.
\textbf{Left}: Model voltage response (black) to variations in $I_{\mathrm{App}}$ (blue). After removal of the inhibitory current, the voltage depolarizes slowly toward its resting potential, exhibits fast oscillations, and then repolarizes gradually over nearly $350$ ms.
\textbf{Right}: Experimental response of cell voltage (black) to an external stimulus $I_{\mathrm{App}}$ (blue) under the PIR protocol. Slow voltage depolarization is followed by a burst of oscillations, and gradual repolarization is observed after the removal of an inhibitory current step. Experimental data from Takakusaki and Kitai~\cite{takakusaki1997ionic} are reproduced with permission from Elsevier.}
\end{figure} 

In the experimental protocol (\Cref{fig:cholt-pir}B), a 300 ms inhibitory current was applied during \textit{in vitro} recordings from a cholinergic PPN neuron with both A-current and LTS~\cite{takakusaki1997ionic}. After the inhibitory step ended, the membrane potential depolarized gradually and produced a burst of oscillations with decreasing amplitude before returning to rest near $-60$ mV gradually over hundreds of milliseconds. The CT model reproduces this qualitative behavior (\Cref{fig:cholt-pir}A), where, after inhibition, the voltage exhibits a delayed depolarization, generates a series of spikes with progressively decreasing amplitude, and then decays back toward baseline. This agreement motivates a closer examination of the characteristic timescales of the low-threshold gating mechanisms of $I_{\mathrm{CaT}}$ and $I_{\mathrm{A}}$, which together shape the delayed and transient nature of the rebound response.

Our analysis for this post-stimulus response will leverage cellular mechanisms on multiple timescales, as introduced in~\Cref{sec:sub:overview-nondim}. The time constants in~\Cref{TAB:ct-nondim} identified by the nondimensionalization procedure reveal three timescales: fast $x = [v_{\mathrm{CT}},m_{\mathrm{Na}},m_{\mathrm{K}},m_{\mathrm{A}},\allowbreak h_{\mathrm{Na}},\allowbreak m_{\mathrm{CaPQ}}]^{\mathsf{T}}$, slow $y = [m_{\mathrm{CaT}},h_{\mathrm{A}}]^{\mathsf{T}}$, and superslow $z=[h_{\mathrm{CaT}},\mathrm{ca}]^{\mathsf{T}}$, and we obtain a nondimensionalized model system of the form~\eqref{sys:ct-nondim}, with $0 < \varepsilon_1, \varepsilon_2 \ll 1$. Further details of this procedure are described in~\Cref{sec:appendix:ct-nondim}. Here, we will introduce definitions of different subsystems and invariant sets in the three-timescale setting used in our analysis of the PIR response of the CT model. 

If we take $\varepsilon_1\rightarrow 0$ in Eq.~\eqref{sys:ct-nondim}, then we obtain the dynamics of the \textit{fast subsystem} in the singular limit, 
\begin{equation}\label{sys:ct-pir-fastsub}
\begin{aligned}
    \frac{dx}{d\tau} &= f(x,y,z), \\
    \frac{dy}{d\tau} &= 0, \\
    \frac{dz}{d\tau} &= 0,
\end{aligned}
\end{equation} 
\noindent in which $y$ and $z$ remain constant. Let the critical manifolds $\mathcal{M}_1^b$ and $\mathcal{M}_1^i$ be defined analogously as in Eq.~\eqref{EQN:M1m} for a three-timescale system with $\bar{I}_{\mathrm{App}}^b=-0.95/(k_v \cdot g_{\mathrm{max}})$ and $\bar{I}_{\mathrm{App}}^i=-3.45/(k_v \cdot g_{\mathrm{max}})$, respectively. Below, we will use the notation $\zeta = (y,z)$ and $\hat{\zeta}$ to denote a frozen subset of components of $\zeta$.

With $\varepsilon_1>0$ in System~\eqref{sys:ct-nondim}, we can take $\varepsilon_2 \rightarrow 0$. We are left with a system where only the right-hand side of $dz/d\tau$ is zero. In this singular limit, only the fast and slow subsystems contribute to the dynamics. The critical points of this reduced system form the \textit{superslow manifold}, defined as 
\begin{equation}
\mathcal{M}_2^m=\{(x,y,z): \; f(x,y,z)=g(x,y)=0, \bar{I}_{\mathrm{App}}=\bar{I}_{\mathrm{App}}^m\}, m\in \{i,e,b\}.
\label{EQN:M2m}
\end{equation}
\noindent In some related instances, we will use the notation $\xi = (x,y)$. 

\begin{figure}[H]
\centering
\includegraphics{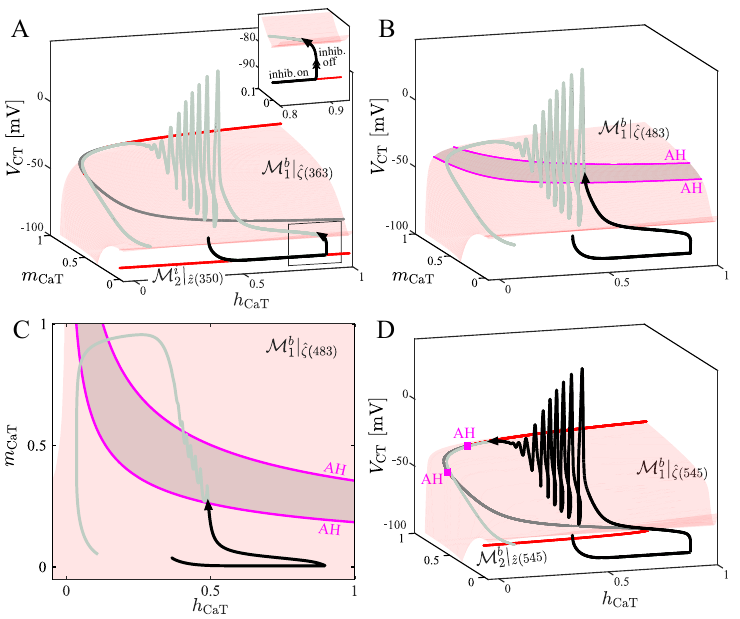}
\caption{\label{fig:cholt-pir-analysis} \small CT model post-inhibitory rebound response is shaped by three-timescale dynamics. Trajectory traces shown in black until $t=363$ (A), $483$ (B,C) and $545$ (D), and shown in gray-green after respective time.
\textbf{A}: Voltage trace from~\Cref{fig:cholt-pir} shown in $(h_{\mathrm{CaT}}, m_{\mathrm{CaT}}, V_{\mathrm{CT}})$-space. $\mathcal{M}_1^b\vert_{\hat{\zeta}(363)}$ is shown at baseline for $\hat{\zeta}(t)=[h_{\mathrm{A}}(t),\mathrm{ca}(t)]$. The superslow manifold $\mathcal{M}_2^i\vert_{\hat{z}(350)}$ comprises two branches of stable (red) and one branch of unstable (gray) equilibria, defined with $\hat{z}(t)=[\mathrm{ca}(t)]$. The inset shows a jump from $\mathcal{M}_2^i$ to $\mathcal{M}_1^b$ in the fast subsystem upon removal of inhibition. The black triangle marks the system state along the trajectory at $t=363$ ms.
\textbf{B}: Onset of oscillations occurs when the trajectory passes a subcritical AH curve (pink) of the fast subsystem that forms a boundary between stable (pink) and unstable (gray) subregions of $\mathcal{M}_1$. $\mathcal{M}_1^b\vert_{\hat{\zeta}(483)}$ is shown and the black triangle indicates system state along the trajectory at $t=483$ ms.
\textbf{C}: Fast variables follow stable periodic orbits before reaching a second branch of supercritical AH bifurcations on $\mathcal{M}_1^b$, after which oscillations conclude. The black triangle indicates the system state at $t=483$ ms.
\textbf{D}: Trajectory reaches superslow manifold $\mathcal{M}_2^b$ (red-stable, gray-unstable). The blue triangle indicates the state in projection at $t=545$ ms. Drift along $\mathcal{M}_2^b$ is governed by the superslow variables, which push the fast-slow subsystem through a supercritical AH and then a subcritical AH bifurcation (pink squares).
}
\end{figure}

In the PIR protocol, after an initial holding period with subthreshold baseline current $\bar{I}_{\mathrm{App}}^b$, inhibition $\bar{I}_{\mathrm{App}}^i$ is applied to the CT neuron. The model trajectory converges to $\mathcal{M}_2^i$ and approaches a stable full system equilibrium at hyperpolarized voltage on this manifold (\Cref{fig:cholt-pir-analysis}A). We define the projection $\mathcal{M}_2^m \vert_{\hat{z}(t)}$ analogously to Eq.~\eqref{EQN:M1myhat}. Let $\hat{z}(t)=[\mathrm{ca}(t)]$ and $\hat{\zeta}(t)=[h_{\mathrm{A}}(t),\mathrm{ca}(t)]$. While $\mathcal{M}_2^m \subset \mathcal{M}_1^m$, notice that $\mathcal{M}_2^m \vert_{\hat{z}(t)} \subset \mathcal{M}_1^m \vert_{\hat{\zeta}(t)}$ only if the slow variable $h_{\mathrm{A}}$ is at its steady-state, satisfying $h_{\mathrm{A}}(t) = h_{\mathrm{A},\infty}(v_{\mathrm{CT}}(t))$. Upon removal of inhibition, dynamics on the fast timescale evolves the trajectory 
from $\mathcal{M}_2^i\vert_{\hat{z}(350)}$ to $\mathcal{M}_1^b\vert_{\hat{\zeta}(363)}$ (\Cref{fig:cholt-pir-analysis}A, inset).

We now introduce \textit{slow time} $\tau_1 = \varepsilon_1 \tau$, where $\tau$ is the fast time. We write the slow time system as

\begin{equation}\label{sys:ct-pir-slowtime}
\begin{aligned}
    \varepsilon_1 \frac{dx}{d\tau_1} &= f(x,y,z), \\
    \frac{dy}{d\tau_1} &= g(x,y), \\
    \frac{dz}{d\tau_1} &= \varepsilon_2 h(x,y,z).
\end{aligned}
\end{equation}

\noindent The dynamics evolving on the slow timescale is considered to be the \textit{slow flow}. Note that on the critical manifold where $f=0$, if we take the singular limit $\varepsilon_2 \rightarrow 0$, then the slow variables $y$ govern the dynamics of the model system.

Once the fast variables jump to a stable branch of equilibria on $\mathcal{M}_1^b$, the slow flow takes the fast subsystem across an AH curve where $\mathcal{M}_1^b$ becomes unstable and the fast variables approach and follow a family of stable periodic orbits (\Cref{fig:cholt-pir-analysis}B,C). At this stage, the slow flow is derived not by assuming $f=0$ but by averaging over the fast oscillations \cite{bertram2017}. The resulting averaged slow flow guides the system toward an attractor of the slow variables. During this excursion, the trajectory continues along the fast subsystem periodic orbits until it crosses a second branch of fast subsystem AH bifurcations on $\mathcal{M}_1^b$, abolishing the oscillations (\Cref{fig:cholt-pir-analysis}B,C). Finally, the slow flow carries the fast variables along $\mathcal{M}_1^b$ until the system reaches the superslow manifold $\mathcal{M}_2^b$ (\Cref{fig:cholt-pir-analysis}D), at which time $dy/d\tau_1$ is near zero. 

Upon entering a small neighborhood of $\mathcal{M}_2^b$, the system transitions to the \textit{superslow flow} on the timescale of $\tau_2 = \varepsilon_2 \tau_1$. The rescaled system

\begin{equation}
\begin{aligned}
    \varepsilon_1 \varepsilon_2 \frac{dx}{d\tau_2} &= f(x,y,z), \\
    \varepsilon_2 \frac{dy}{d\tau_2} &= g(x,y), \\
    \frac{dz}{d\tau_2} &= h(x,y,z)
\end{aligned}\label{sys:ct-pir-superslowtime}
\end{equation}
\noindent exposes the superslow flow for a system on the superslow manifold $\mathcal{M}_2^b$ where $f=g=0$. On $\mathcal{M}_2^b$ in the singular limit, the drift in superslow variables now approaches an attractor of the superslow subsystem (\Cref{fig:cholt-pir-analysis}D). The trajectory, governed by the superslow variables, passes through a supercritical AH bifurcation in the fast-slow subsystem, causing a delayed bifurcation effect \cite{desroches2012} in which the fast-slow variables remain near a family of unstable equilibria along $\mathcal{M}_2^b$, without an instant onset of oscillations. Before oscillations have a chance to begin, the superslow flow takes the system through a second, subcritical fast-slow AH bifurcation, after which the trajectory drifts along the stable branch of $\mathcal{M}_2^b$ before reaching the attractor of the superslow subsystem.

In summary, this three-timescale decomposition demonstrates how the rebound oscillations emerge and disappear from the interplay of three distinct timescales. The slow dynamics guides the fast variables through the onset, persistence and offset of oscillations, while the superslow dynamics subsequently steers the system through a delayed AH bifurcation followed by another AH bifurcation without causing the onset of oscillations. The superslow drift then helps shape the gradual repolarization before reaching a full-system equilibrium state, at which time the PIR response is concluded. At each stage, we are able to offer a mechanistic explanation for the experimentally observed features of the PIR activity.

\vspace{0.1cm}
\textit{\textbf{Remark}} According to System~\eqref{sys:ct-pir-slowtime}, on the critical manifold $\mathcal{M}_1^b$, the slow subsystem should govern the dynamics until the slow variables reach a small neighborhood of the superslow manifold $\mathcal{M}_2^b$. However, between the fast jump at the removal of inhibition (\Cref{fig:cholt-pir}A) and the crossing of the fast subsystem AH curve (\Cref{fig:cholt-pir-analysis}B), we see a more significant drift in superslow variable $h_{\mathrm{CaT}}$, decreasing from 0.9 to 0.5, than in slow variable $m_{\mathrm{CaT}}$. The timescale distinction between slow and superslow variables breaks down at low voltages because the nonlinear dependence of $\tau_{m_{\mathrm{CaT}}}(v_{\mathrm{CT}}),\tau_{h_{\mathrm{CaT}}}(v_{\mathrm{CT}})$ and $m_{\mathrm{CaT},\infty}(v_{\mathrm{CT}}), h_{\mathrm{CaT},\infty}(v_{\mathrm{CT}})$ on $v_{\mathrm{CT}}$ cannot be fully captured by the time constants $R_{m_{\mathrm{CaT}}},R_{h_{\mathrm{CaT}}}$ chosen to identify the subsystems evolving on distinct timescales. We do not go into further details of this blurring of timescales here, since these would add complication without clarifying the main mechanisms underlying the PIR response.
\vspace{0.1cm}

\subsection{NC and CT models predict post-inhibitory facilitation in PPN cells with low-threshold calcium channels}\label{sec:nc-ct-pif}

Beyond reproducing known PPN dynamics and suggesting the underlying mechanisms, our models also allow for the exploration of new stimulation protocols to predict how the different PPN cell types may respond. 
Here, we consider a stimulus protocol associated with an effect called post-inhibitory facilitation (PIF). PIF is a phenomenon observed in neurons with Type III excitability, in which a brief inhibitory input enhances the effectiveness of a subsequent excitatory stimulus \cite{prescott2008biophysical,meng2012type,rinzel2013nonlinear,huguet2017phasic}. Rather than triggering rebound spiking directly, the inhibitory pulse transiently recruits subthreshold currents whose slow gating dynamics create a window of increased excitability. Motivated by this mechanism, we apply an inhibitory–excitatory pulse protocol to our PPN models to predict how different intrinsic conductances will shape stimulus integration. We will show that both models containing low-threshold calcium current (NC and CT) exhibit clear PIF, and we will explain how their intrinsic currents shape these responses.

\begin{figure}[H]
\centering
\includegraphics{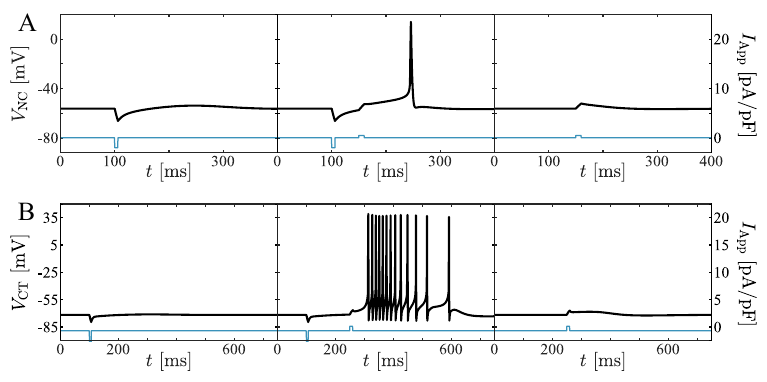}
\caption{\label{fig:nonchol-cholt-pif} \small Post-inhibitory facilitation responses in NC and CT PPN models.
\textbf{A}: NC model exhibits post-inhibitory facilitation. The voltage (black) starts near a resting potential of $-60$ mV with $I_{\mathrm{App}}=0$ pA/pF. A brief inhibitory stimulus applied at $t=100$ ms followed by a small excitatory stimulus at $t=150$ ms induces a single spike in the model voltage (center). In comparison, the inhibitory pulse alone (left) or excitatory pulse alone (right) does not elicit spiking in the model cell.
\textbf{B}: CT model exhibits a post-inhibitory spiking response. The voltage (black) starts near a resting potential near $-60$ mV with $I_{\mathrm{App}}=0$ pA/pF. A brief inhibitory stimulus applied at $t=100$ ms followed by a small excitatory stimulus at $t=250$ ms induces a series of spikes with decreasing frequency before returning to rest (center). In comparison, having only the inhibitory pulse (left) or only the excitatory pulse (right) does not elicit spiking in the model cell.}
\end{figure}

Starting from the stable rest state in the NC model, we apply a brief inhibitory pulse followed, after a short delay, by a brief excitatory pulse. Under this protocol, the model produces a single spike~(\Cref{fig:nonchol-cholt-pif}A, center), even though neither the inhibitory pulse alone~(\Cref{fig:nonchol-cholt-pif}A, left) nor the excitatory pulse alone~(\Cref{fig:nonchol-cholt-pif}A, right) is sufficient to elicit spiking. Applying the same protocol to the CT PPN model also produces PIF, but in this case the response consists of 13 spikes over 205 ms ~(\Cref{fig:nonchol-cholt-pif}B, center). As in the NC model, neither the inhibitory pulse alone~(\Cref{fig:nonchol-cholt-pif}B, left) nor the excitatory pulse alone~(\Cref{fig:nonchol-cholt-pif}B, right) evokes spiking activity. Applying the inhibitory–excitatory pulse protocol to the C model did not elicit any spikes, even when pulse durations and inter-stimulus delays were varied (not shown). This outcome is consistent with the known properties of the C model, as its only low-threshold current is the potassium current $I_{\mathrm{A}}$. In this case, the inhibitory pulse hyperpolarizes the cell and transiently activates $I_{\mathrm{A}}$. The resulting period of increased $I_{\mathrm{A}}$ conductance suppresses excitability, preventing subsequent excitatory pulses from inducing spiking.

\subsubsection{\textbf{Three-timescale switching in NC PIF response}}\label{sec:nc-pif}

To analyze PIF in the NC model, we follow the nondimensionalization procedure as outlined in~\Cref{sec:sub:overview-nondim,sec:appendix:ct-nondim}.
To capture the PIF dynamics in terms of a timescale decomposition, it turns out that we need to split the trajectory into an initial, subthreshold phase and a later, spiking phase. The time constants for the model variables in these phases appear as ``timescale-1'' and ``timescale-2'', respectively, in~\Cref{tab:nc-timeconstants}.
To represent the shift in the timescales and the roles of $h_{\mathrm{Na}}$ and $m_{\mathrm{K}}$ across these phases, we start from the classification used for the NC model in~\Cref{sec:nc-sdp} except that we split off $h_{\mathrm{Na}}$ and $m_{\mathrm{K}}$ into their own group.
That is, we represent the model as a three-timescale system $(x,y,z)$ with fast $x=[v_{\mathrm{NC}}, m_{\mathrm{Na}},m_{\mathrm{CaPQ}}]^{\mathsf{T}}$, slow $y = [h_{\mathrm{Na}},m_{\mathrm{K}}]^{\mathsf{T}}$, and superslow variables $z=[m_{\mathrm{CaT}}, h_{\mathrm{CaT}}, \mathrm{ca}]^{\mathsf{T}}$.
With this decomposition, our system can be written in the form~\eqref{sys:ct-nondim}. As we proceed with the analysis, we will group the slow variables $y$ with either $x$ or $z$, as appropriate.

\begin{figure}[H]
\centering
\includegraphics{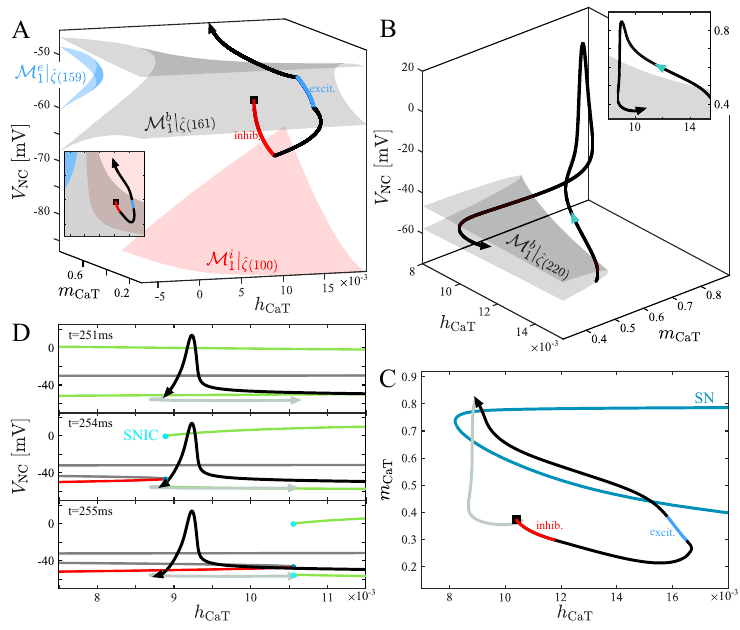}
\caption{\label{fig:nonchol-pif-analysis} \small Analysis of PIF response in NC model leverages three-timescale dynamics.
\textbf{A}: NC model trajectory during inhibition (red line), during excitation (blue line), and at baseline (black line). Projections of the critical manifolds, $\mathcal{M}_1^i\vert_{\hat{\zeta}(100)}$ (pink surface), $\mathcal{M}_1^e\vert_{\hat{\zeta}(159)}$ (blue surface) and $\mathcal{M}_1^b\vert_{\hat{\zeta}(161)}$ (gray surface) for $\hat{\zeta}(t)=[m_{\mathrm{K}}(t),h_{\mathrm{Na}}(t),\mathrm{ca}(t)]$, of the $(x,\zeta)$ decomposition for subthreshold activity are shown. The black square indicates the point along the trajectory at the onset of inhibition.
The inset shows projection down onto the $(h_{\mathrm{CaT}},m_{\mathrm{CaT}})$-plane with the same axes and tick marks.
\textbf{B}: $\mathcal{M}_1^b$ (gray) is shown for $\hat{\zeta}(220)$, the time marked by the teal triangle on the system trajectory (black), which is shown between $t=160$ and $500$ ms. The trajectory remains close to the fold of $\mathcal{M}_1^b$ during slow depolarization. The inset shows the projection onto the $(h_{\mathrm{CaT}},m_{\mathrm{CaT}})$-plane with the same axis units.
\textbf{C}: Slow depolarization along the fast subsystem SN curve (teal) that occurs before the large-amplitude oscillation in dimensionalized $V_{\mathrm{NC}}$. The teal curve represents the projection of $\mathcal{B}^{\mathrm{SN}}$ onto the dimensionalized phase space. The trajectory, starting at $t=0$ ms (black square), is also shown for $t<249$ (red, black, blue following color scheme in A) along with its later path (gray-green). The branch of SN is computed for values of the slow variables between $t=160$ and $249$ ms.
\textbf{D}: Large-amplitude oscillation terminates at SNIC crossing of fast variables $\xi$ in $(\xi,z)$ decomposition. The red/gray curves form the critical manifold of the two-timescale $(\xi,z)$ system; a family of periodic orbits of this system is shown in green. The trajectory completes one oscillation along this periodic orbit family and then crosses a SNIC bifurcation (center panel) of the $\xi$ dynamics such that no additional spikes occur. The $\xi$ components then approach and follow the stable branch of $\mathcal{M}_2^b$ (bottom panel), where the slow flow takes over.
}
\end{figure}

Three levels of current are applied during the PIF protocol. At baseline, after nondimensionalization (\Cref{sec:appendix:ct-nondim}), $\bar{I}_{\mathrm{App}}^b=0$, leading to the definition of the critical manifold $\mathcal{M}_1^b$ by Eq.~\eqref{EQN:M1m}. During inhibition and excitation, $\bar{I}_{\mathrm{App}}^i=-2/(k_v \cdot g_{\mathrm{max}})$ and $\bar{I}_{\mathrm{App}}^e=0.5/(k_v \cdot g_{\mathrm{max}})$, respectively, and the critical manifolds are $\mathcal{M}_1^i$ and $\mathcal{M}_1^e$. We define the superslow manifold $\mathcal{M}_2^b$ using Eq.~\eqref{EQN:M2m} with applied current valued at baseline.

For subthreshold activity with $v_{\mathrm{NC}} \in [-0.66,-0.45]$, we introduce the timescale grouping $(x,(y,z))$ and again use $\zeta=(y,z)$.
Before either of the stimulation pulses, the model trajectory lies at an equilibrium near $v_{\mathrm{NC}}=-0.53$. Once inhibition is applied to the model cell, the fast variables are attracted to a low-voltage branch of $\mathcal{M}_1^i$ instead, near $v_{\mathrm{NC}}=-0.8$. When the inhibition is removed, the fast variables approach the stable branch of $\mathcal{M}_1^b$ near $-0.58$. The subsequent excitatory pulse yields a stable branch of $\mathcal{M}_1^e$ at higher voltages near $v_{\mathrm{NC}}=-0.5$, which also attracts the fast variables. Before the fast variables can make much progress toward $\mathcal{M}_1^e$, the excitation is removed, and the fast variables once again approach $\mathcal{M}_1^b$. These steps are illustrated in~\Cref{fig:nonchol-pif-analysis}A. 

After the two current pulses, the model voltage exhibits a gradual depolarization that eventually leads into a single spike. On $\mathcal{M}_1^b$, $f_1(x,y,z)=0$, where $f_1$ denotes the right-hand side of the $v_{\mathrm{NC}}$ equation in System~\eqref{sys:nc-ramp-fasttime} (see~\Cref{sec:appendix-nc-nondim}). Multiple timescale analysis suggests that, on $\mathcal{M}_1^b$, the slow flow of System~\eqref{sys:ct-pir-slowtime} governs trajectory behavior. As shown in~\Cref{fig:nonchol-pif-analysis}B, this dynamics pushes the fast variables past the fold on $\mathcal{M}_1^b$, yet $v_{\mathrm{NC}}$, a fast variable, does not immediately jump to another attractor. Instead, we see that the dynamics follows the fold branch, where $v_{\mathrm{NC}}$ appears to be pinned near $f_1=0$. \Cref{fig:nonchol-pif-analysis}C illustrates the slow drift along
$$\mathcal{B}^{\mathrm{SN}} = \bigcup_{\tau\in \tau_{\mathrm{range}}} \left\{(x,\zeta) \in \mathcal{M}_1^b \; \vert \; \frac{df(x,\zeta(\tau))}{dv} = 0 \right\}$$

\noindent for $\tau_{\mathrm{range}}=[0.16,0.249]$. The three-timescale decomposition can break down when the system is tangent to a curve of knees in the fast subsystem. As a result, fast variable $v_{\mathrm{NC}} \in x$ evolves more slowly than its timescale would suggest, closer to the timescales of $h_{\mathrm{Na}}$ and $m_{\mathrm{K}}$.

\vspace{0.1cm}
\textit{\textbf{Remark}}: We can disrupt this balance by increasing the timescale of $m_{\mathrm{CaT}}$ to $t_{m_{\mathrm{CaT}}}=5$ ms. With this adjustment, we no longer find the drift along the curve of knees, and the trajectory quickly moves away from the critical manifold under the fast flow. In addition, setting $m_{\mathrm{Na}}=m_{\mathrm{Na},\infty}(v_{\mathrm{NC}})$ and $m_{\mathrm{CaPQ}}=m_{\mathrm{CaPQ},\infty}(v_{\mathrm{NC}})$ for all $v_{\mathrm{NC}}$ does not qualitatively alter the trajectory, which confirms that these other fast variables are not part of this mixing of timescales.
\vspace{0.1cm}

As the voltage depolarizes above $v_{\mathrm{NC}}=-0.45$, the $(x,\zeta)$ decomposition loses validity. For $v_{\mathrm{NC}} \in [-0.45,0.14]$, slow variables $h_{\mathrm{Na}}, m_{\mathrm{K}}$ evolve at a faster rate (see \Cref{tab:nc-timeconstants}). As a result, we view the three-timescale system from the perspective of $((x,y),z)$ for the remainder of the analysis of this response, and introduce $\xi=(x,y)$. The trajectory approaches and follows a stable family of periodic orbits of the $\xi$ dynamics for one period before dynamics in the superslow variables, $z$, takes the fast variables across the SNPO bifurcation in $\mathcal{M}_2^b$ that terminates the oscillations (\Cref{fig:nonchol-pif-analysis}D). Variables in $\xi$ then approach the nearby attracting branch of $\mathcal{M}_2^b$ around $v_{\mathrm{NC}}=-0.6$. After the $\xi$ flow reaches this stable branch of equilibria, the superslow dynamics takes over as the $z$ components approach their attractor on $\mathcal{M}_2^b$. This final superslow drift returns the system to its resting potential.

\subsubsection{\textbf{PIF oscillations are based on three-timescale dynamics in CT model}}\label{sec:ct-pif}

For the full CT system described in System~\eqref{SYS:PPN} with Eq.~\eqref{EQN:CHOLT}, after the nondimensionalization steps analogous to the procedure described in~\Cref{sec:sub:overview-nondim,sec:appendix:ct-nondim}, we identify (see~\Cref{TAB:ct-nondim}) 

\begin{figure}[ht!]
\centering
\includegraphics{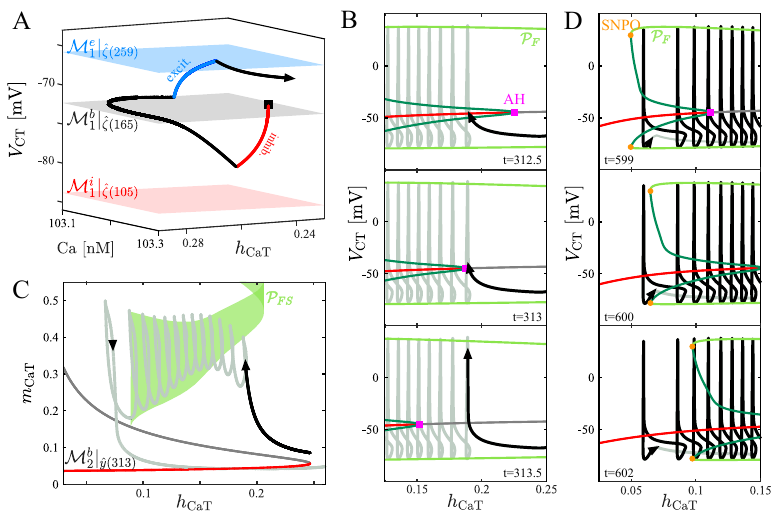}
\caption{\label{fig:cholt-pif-analysis} PIF response in CT model is explained through three-timescale decomposition.
\textbf{A}: Critical manifolds $\mathcal{M}_1^i\vert_{\hat{\zeta}(105)}$, $\mathcal{M}_1^e\vert_{\hat{\zeta}(259)}$, and $\mathcal{M}_1^b\vert_{\hat{\zeta}(165)}$ shown during inhibitory pulse (red), during excitatory pulse (blue) and at baseline (black), respectively, with $\hat{\zeta}(t)=[m_{\mathrm{CaT}}(t),h_{\mathrm{A}}(t)]$. The black square indicates the point along the trajectory at the onset of inhibition.
\textbf{B}: During onset of oscillations, the trajectory crosses a fast subsystem AH bifurcation (pink squares corresponding to AH locations at different points in time) on $\mathcal{M}_1^b$. In each panel, the trajectory up until the indicated time is drawn in black and the future path of the trajectory is shown in gray-green. Stable periodic orbits $\mathcal{P}_{F}$ (min/max in light green) are drawn for the fast subsystem, with those variables not shown fixed at their values at the time $t$ indicated in each panel. Red and gray segments on $\mathcal{M}_1^b$ indicate stable and unstable equilibria, respectively. AH bifurcation in fast subsystem at $t=313$ ms precedes the arrival of the fast-slow subsystem at its family of stable periodic orbits in panel C.
\textbf{C}: Projection of trajectory for $t\in(265,313)$ shown in black, and for $t \in (313,900)$ shown in gray-green. The black triangles near $h_{\mathrm{CaT}}=0.2, 0.08$ indicate points on trajectory at $t=313, 600$ ms, respectively. $\mathcal{M}_2^b\vert_{\hat{z}(313)}$ represents the superslow manifold where $\hat{z}(t)=[\mathrm{ca}(t)]$. The family of periodic orbits $\mathcal{P}_{FS}$ of the fast-slow subsystem is shown in green. The SNPO crossing of the fast-slow subsystem occurs prior to the fast subsystem crossing of the SNPO at $t=600$ ms, shown in panel D.
\textbf{D}: During the offset of oscillations, the fast subsystem crosses its SNPO bifurcation, and stable periodic orbits are no longer present at low values of $h_{\mathrm{CaT}}$. The color coding of the invariant sets follows the same schematic as in B.
}
\end{figure}

\noindent $x=[v_{\mathrm{CT}},m_{\mathrm{Na}}, m_{\mathrm{K}},m_{\mathrm{A}}, \allowbreak h_{\mathrm{Na}}, m_{\mathrm{CaPQ}}]^{\mathsf{T}}$ as the fast subsystem, while the remaining variables operate on slower timescales and form the slow $y=[m_{\mathrm{CaT}},h_{\mathrm{A}}]^{\mathsf{T}}$ and superslow subsystems $z=[h_{\mathrm{CaT}},\mathrm{ca}]^{\mathsf{T}}$, exactly as in our analysis of post-inhibitory oscillations in the CT model in~\Cref{sec:ct-pir}. This three-timescale decomposition yields a system of the form~\eqref{sys:ct-nondim}. As in the previous subsection, we will see that it is appropriate to group certain timescale components together on segments of the trajectory of interest that span distinct voltage ranges.

Before the application of stimuli, the trajectory lies near a critical point with $v_{\mathrm{CT}} \approx -0.72$.
During the inhibitory stimulus, the critical manifold $\mathcal{M}_1^i$ is defined for an applied current of $\bar{I}_{\mathrm{App}}^i=-2.7/(k_v \cdot g_{\mathrm{max}})$. The fast variables approach $\mathcal{M}_1^i$ on the fast, $\tau$ timescale until the inhibition is removed (\Cref{fig:cholt-pif-analysis}A). Then, the fast variables approach $\mathcal{M}_1^b$, and the slow drift briefly takes over, on the slow timescale $\tau_1$, until the excitation is applied. The attractor of the fast variables next becomes $\mathcal{M}_1^e$, defined with $\bar{I}_{\mathrm{App}}^e=0.15/(k_v \cdot g_{\mathrm{max}})$, toward which the fast subsystem evolves until excitation is removed and the fast variables again approach $\mathcal{M}_1^b$.

In the singular limit, once the fast variables reach $\mathcal{M}_1^b$, the drift in the slow variables again takes over, and it pulls the trajectory through a fast subsystem AH curve on the critical manifold (\Cref{fig:cholt-pif-analysis}B) near $t=313$ ms. At this time, the fast variables begin to oscillate as they follow a family of stable periodic orbits. While the fast subsystem has this family of stable periodic orbits, call it $\mathcal{P}_F$, the larger fast-slow subsystem also has a distinct such family, call it $\mathcal{P}_{FS}$ (\Cref{fig:cholt-pif-analysis}C, green). In the next phase of the dynamics, the fast variables oscillate along $\mathcal{P}_F$ while the slow flow, averaged around these oscillations, carries the trajectory toward $\mathcal{P}_{FS}$. Once convergence to $\mathcal{P}_{FS}$ has occurred, the superslow variables take over the dynamics, with drift on the superslow timescale $\tau_2$ (see System~\eqref{sys:ct-pir-superslowtime}), while the fast-slow variables follow $\mathcal{P}_{FS}$. 

Near $h_{\mathrm{CaT}}=0.9$, the trajectory crosses a SNIC bifurcation of the 8-dimensional fast-slow subsystem, terminating in a homoclinic orbit. The fast-slow subsystem no longer exhibits stable periodic orbits, and we can see an increase in period with each oscillation (\Cref{fig:nonchol-cholt-pif}B, center). Yet the fast subsystem continues to oscillate as the trajectory follows the family of fast subsystem periodic orbits $\mathcal{P}_{F}$ (\Cref{fig:cholt-pif-analysis}D, top). The drift along this family is governed by the slow variables, with dynamics based on averaging over the fast periodic orbits, and pulls the trajectory through a fast subsystem SNPO bifurcation (\Cref{fig:cholt-pif-analysis}D, center) near time $t=600$ ms. This point marks the termination of oscillations, after which the fast variables approach a stable branch of the critical manifold (\Cref{fig:cholt-pif-analysis}D, bottom) and the slow variables subsequently settle to their attractor on $\mathcal{M}_2^b$ (\Cref{fig:cholt-pif-analysis}C).

To conclude, we note that certain phases of the CT PIF response can be characterized using the attractors of the combined fast-slow subsystem and superslow subsystem. However, we see that the fast variables reach their attracting periodic orbits before the slow variables reach their attractor. This short delay calls for the distinction between variables in $x$ and $y$. Similarly, the fast-slow subsystem crosses a SNIC bifurcation before the stable oscillations in the fast subsystem terminate. The offset of oscillations is thus another example of dynamics that would not be captured by a two-timescale decomposition of this model system. In~\Cref{fig:cholt-pif-analysis}B-D, we can also see oscillations of non-negligible magnitude in the fast $v_{\mathrm{CT}}$ and slow $m_{\mathrm{CaT}}$ variable, whereas with superslow variable $h_{\mathrm{CaT}}$, there is an approximately monotone drift with little change in magnitude during each spike in $v_{\mathrm{CT}}$.

\section{Discussion}
\label{sec:discussion}

\noindent In this work, we introduce and analyze three distinct models of neurons in the PPN, based on transmembrane voltage dynamics observed experimentally in response to electrical stimulation patterns. This is the first set of models for multiple cell classes in PPN, calibrated to a range of data specific to PPN. In each model, we show that the inclusion of specific ion currents allows for various experimental observations to be reproduced, using a single set of currents and parameters for all simulations. The NC model includes high- and low-voltage activated calcium channels, in addition to calcium-activated potassium channels, to successfully capture small-voltage gamma oscillations in response to a ramping excitatory current. Transient spiking activity during a small excitatory step current is reproduced by the same set of currents in the NC model as well. The C model, after inhibition, exhibits a delay before firing, characteristic of PPN cells with a low-voltage activated potassium A-current. Lastly, the CT model expresses low-voltage activated potassium and calcium currents that yield rebound spiking activity when subjected to an inhibitory current step (PIR) or sequential inhibitory and excitatory current steps (PIF). In our analysis of the PIF rebound activity, oscillations in the fast subsystem terminate through an SNPO bifurcation, whereas in the combined fast-slow subsystem the periodic orbits instead end in a SNIC bifurcation. In the fast subsystem, the A-current activation variable $m_{\mathrm{A}}$ is treated as dynamic, while in the fast-slow system, both activation $m_{\mathrm{A}}$ and inactivation $h_{\mathrm{A}}$ of the A-current evolve. Slow outward currents such as the A-current are known to alter firing onset~\cite{connor1971prediction} and can shift the onset mechanisms to a SNIC bifurcation~\cite{koch1998methods}.

In addition to revealing key mechanisms underlying each of these voltage responses, our models can provide specific predictions about PPN cell responses to new current protocols or about activity expected to arise if specific currents are blocked. For example, we showed that our two models with T-type calcium channels (CT and NC) exhibit distinct post-inhibitory facilitation (PIF) effects, where an inhibitory-excitatory pulse pair evokes transient activity, in the absence of any large-amplitude voltage response to either pulse alone. These dynamics represent predictions about expected behaviors of CT and NC PPN neurons. Interestingly, although PIF is typically associated with neurons with Type III excitability \cite{prescott2008biophysical,meng2012type,rinzel2013nonlinear}, recent mathematical work has shown that it is actually a more general phenomenon \cite{rubin2021}, and our findings provide an example of this generalization.

We show how each of these qualitative spike patterns emerges from multiple timescale dynamics. Typical multiple timescale analysis proceeds via the identification or designation of two, or in some cases three, distinct timescales that shape the dynamics exhibited by a nonlinear system. Our analyses, in contrast, highlight the fact that some variables may act on different timescales on distinct regions of phase space. Correspondingly, we find that fluidity in the designation of timescale subsystems is necessary to explain the observed PIF dynamics \cite{john2024novel}. Such shifting in the relative timescales of model variables arises as a natural consequence of temporal scaling functions that depend on one or more modeled quantities, in this case voltage, in biological systems. Hence, this timescale fluidity may represent a widespread feature of neural models, for which a more systematic mathematical treatment would be a useful direction for future work. 
Another issue for future analysis is whether there are clear conditions for how variables should be split into distinct timescale classes. Some of the groupings that we find provide an accurate representation of model dynamics in terms of timescale decomposition and involve splitting variables with similar time constants into distinct classes (e.g., see~\Cref{tab:nc-timeconstants}), perhaps reflecting the fact that some variables evolve at different rates along different parts of relevant solutions.

The three models presented here are the first single-cell PPN models to be calibrated to reproduce responses across a range of experimental current protocols. Due to limitations in the scope of published data, once the current set for each model was selected, individual channel formulations and initial parameter values were drawn from measurements across different animals and brain regions. In addition, although noncholinergic PPN cells include GABAergic and glutamatergic subtypes, limitations in subtype-specific data led us to represent them with a single NC model. Each model uses one parameter set to capture a range of dynamical responses; however, given the high-dimensional parameter space, there may in theory exist other parameter sets that yield similar dynamics in these models. Our bifurcation analyses suggest that the qualitative dynamics of our models are robust to moderate parameters changes, but a systematic exploration of parameter space and full sensitivity analysis are beyond the scope of this study. To circumvent issues of heterogeneity in, and lack of data about, capacitance in PPN cells, we present model conductance values in per-capacitance units (nS/pF). In these units, conductances in the NC model are smaller than those of the cholinergic models, and we find that relative maximal conductance values, more than absolute conductances, primarily govern the observed activity patterns.

Future work could include the study of interactions of PPN neurons with brain regions to which it is coupled synaptically, such as the SNr and STN in the basal ganglia \cite{fallah2024inhibitory,falasconi2025,hu2026multimodal}. These pathways likely contribute to motor functions and may be sites of altered dynamics in disorders involving motor effects, such as Parkinson's disease. Recent experiments have shown that transient optogenetic stimulation, when applied to subpopulations within the basal ganglia, can result in long-lasting motor rescue in dopamine-depleted rodents~\cite{mastro2017,cundiff2024}. 
Given that PPN neurons are a direct target of outputs from basal ganglia neurons in motor-related pathways, our PPN models can be leveraged to analyze the mechanisms behind the motor effects of these stimulation experiments, as a step toward development of new therapeutic approaches for PD patients. Moreover, cholinergic PPN neurons are known to undergo apoptosis in PD, and the future incorporation of our C and CT PPN models into a basal ganglia-PPN circuit model would allow for the study of how the loss of PPN neurons impacts motor function.

\appendix

\section{An introduction to multiple timescale analysis\\ and why we need it}\label{apn:GSPT}

\noindent Biological systems and their components frequently operate across multiple distinct time- scales. In neurons, for example, the processes that shape membrane voltage evolve at different rates. Sodium channel activation and the associated influx of ions occur over just a few milliseconds, producing rapid changes in membrane potential. In contrast, intracellular calcium concentration typically accumulates and decays over hundreds of milliseconds, changing only gradually during fast voltage events. While fast mechanisms can generate nearly instantaneous effects, slower processes tend to modulate activity over longer durations and contribute to sustained changes within the cell. The interaction of processes across these timescales introduces substantial complexity and gives rise to a rich variety of dynamical behavior. In mathematical modeling, multi-timescale phenomena are often studied by splitting variables into those that evolve rapidly (``fast'') and those that evolve more gradually (``slow''). This framework, known as multiple timescale decomposition, provides a systematic way to analyze dynamical systems arising in biology and neuroscience. If we denote the fast variables by $x$ and the slow variables by $y$, we can temporarily exaggerate the difference in their rates by treating the slow variables as effectively constant. Under this approximation, the $y$ variables are considered stationary, allowing us to isolate and analyze the dynamics of $x$, determine their steady states, and examine how these steady states vary under changes in $y$. After understanding the dynamics of the fast variables and identifying the relevant attractor of the fast dynamics, we restore the gradual evolution of the slow variables that occurs with the fast variables on this attractor to study how they drift over time and shape the full system dynamics. This perspective clarifies how rapid events, such as action potentials, can be modulated by slower processes like calcium dynamics or neuromodulatory changes, without requiring all variables to be analyzed simultaneously in a high-dimensional setting.

\section{Nondimensionalization}
\label{sec:appendix-nondim}

\noindent In this Appendix, we outline the steps used in the nondimensionalization of the CT model for the analysis of its post-inhibitory oscillatory activity. We then offer summaries of analogous procedures for timescale identification in the C and NC models. A similar procedure applied in the analysis of a three-timescale system can be found in Nan et al.~\cite{nan2015understanding}. Further examples of the nondimensionalization process are provided in other past works~\cite{bertram2017,kuehn2015multiple,rubin2007giant}.

\subsection{Multiple timescales of the CT PPN model}\label{sec:appendix:ct-nondim}

In~\Cref{sec:sub:overview-nondim}, we provided an overview of the nondimensionalization procedure for the identification of three subsystems evolving on distinct timescales in the CT model, allowing us to rewrite the model as the three-timescale System~\eqref{sys:ct-nondim} using the dimensionless variables $v_{\mathrm{CT}}, c, \tau$. Here, we present the computations needed for this nondimensionalization, including rescaling parameters and computing time constants $R_i$ in~\Cref{TAB:ct-nondim}.

To rescale the dynamics of an arbitrary gating variable $p$ of the CT model, we define $T_{p} = \max_{V_{\mathrm{CT}} \in V_{\mathrm{range}}}\{1/t_p(V_{\mathrm{CT}})\}$ to obtain a rescaled, dimensionless version of the gating variable time constant $\tau_p(v_{\mathrm{CT}}) = T_p t_p(v_{\mathrm{CT}})$ over the relevant voltage range. We define $g_{\mathrm{max}} = \max_{i \in \mathcal{I_{\mathrm{CT}}}} \{g_i\}$ and see that $g_{\mathrm{max}} = g_{\mathrm{Na}} = 50$ nS. Substituting the dimensionless variables and newly defined parameters into System~\eqref{SYS:PPN} for $v_{\mathrm{CT}}$ with $\mathcal{I}_{\mathrm{CT}}$ currents results in the following dimensionless system:

\begin{equation}\label{sys:appendix-ct-nondim1}
\begin{aligned}
    \frac{dv}{d\tau} =& \; R_v \bar{g}_{\mathrm{Na}} \cdot m_{\mathrm{Na}}^3 \cdot h_{\mathrm{Na}} \cdot (v_{\mathrm{CT}} - \bar{E}_{\mathrm{Na}}) \; - \;  
    \bar{g}_{\mathrm{K}} \cdot m_{\mathrm{K}}^4\cdot (v_{\mathrm{CT}} - \bar{E}_{\mathrm{K}})\\
    &\; - \; 
    \bar{g}_{\mathrm{L}}\cdot (v_{\mathrm{CT}} - \bar{E}_{\mathrm{L}}) \; - \; \bar{g}_{\mathrm{CaT}}\cdot m_{\mathrm{CaT}}^2 \cdot h_{\mathrm{CaT}} \cdot (v_{\mathrm{CT}} - \bar{E}_{\mathrm{Ca}})\\
    &\; - \; \bar{g}_{\mathrm{CaPQ}} \cdot m_{\mathrm{CaPQ}} \cdot (v_{\mathrm{CT}} - \bar{E}_{\mathrm{Ca}}) \; - \; \bar{g}_{\mathrm{A}} \cdot m_{\mathrm{A}}^2 \cdot h_{\mathrm{A}} \cdot (v_{\mathrm{CT}} - \bar{E}_{\mathrm{K}})\\
    &\; - \; \bar{g}_{\mathrm{KCa}} \cdot m_{\mathrm{KCa}}(v_{\mathrm{CT}},\mathrm{ca}) \cdot (v_{\mathrm{CT}} - \bar{E}_{\mathrm{K}}) \; - \; \bar{I}_{\mathrm{App}}, \\
    \frac{dm_{\mathrm{Na}}}{d\tau} =& \; R_{m_{\mathrm{Na}}} \frac{m_{\mathrm{Na},\infty}(v_{\mathrm{CT}})-m_{\mathrm{Na}}}{\tau_{m_{\mathrm{Na}}}(v_{\mathrm{CT}})},\\
    \frac{dh_{\mathrm{Na}}}{d\tau} =& \; R_{h_{\mathrm{Na}}} \frac{h_{\mathrm{Na},\infty}(v_{\mathrm{CT}})-h_{\mathrm{Na}}}{\tau_{h_{\mathrm{Na}}}(v_{\mathrm{CT}})},\\
    \frac{dm_{\mathrm{K}}}{d\tau} =& \; R_{m_{\mathrm{K}}} \frac{m_{\mathrm{K},\infty}(v_{\mathrm{CT}})-m_{\mathrm{K}}}{\tau_{m_{\mathrm{K}}}(v_{\mathrm{CT}})},\\
    \frac{dm_{\mathrm{CaT}}}{d\tau} =& \; R_{m_{\mathrm{CaT}}} \frac{m_{\mathrm{CaT},\infty}(v_{\mathrm{CT}})-m_{\mathrm{CaT}}}{\tau_{m_{\mathrm{CaT}}}(v_{\mathrm{CT}})},\\
    \frac{dh_{\mathrm{CaT}}}{d\tau} =& \; R_{h_{\mathrm{CaT}}} \frac{h_{\mathrm{CaT},\infty}(v_{\mathrm{CT}})-h_{\mathrm{CaT}}}{\tau_{h_{\mathrm{CaT}}}(v_{\mathrm{CT}})},\\
    \frac{dm_{\mathrm{CaPQ}}}{d\tau} =& \; R_{m_{\mathrm{CaPQ}}} \frac{m_{\mathrm{CaPQ},\infty}(v_{\mathrm{CT}})-m_{\mathrm{CaPQ}}}{\tau_{m_{\mathrm{CaPQ}}}(v_{\mathrm{CT}})},\\
    \frac{dm_{\mathrm{A}}}{d\tau} =& \; R_{m_{\mathrm{A}}} \frac{m_{\mathrm{A},\infty}(v_{\mathrm{CT}})-m_{\mathrm{A}}}{\tau_{m_{\mathrm{A}}}(v_{\mathrm{CT}})},\\
    \frac{dh_{\mathrm{A}}}{d\tau} =& \; R_{h_{\mathrm{A}}} \frac{h_{\mathrm{A},\infty}(v_{\mathrm{CT}})-h_{\mathrm{A}}}{\tau_{h_{\mathrm{A}}}(v_{\mathrm{CT}})},\\
    \frac{d\mathrm{ca}}{d\tau} =& \; R_{\mathrm{Ca}}
    f_{\mathrm{Ca}} \left( 
    \frac{k_{\tau}}{k_{\mathrm{ca}}} \frac{\bar{g}_{\mathrm{CaT}}\cdot m_{\mathrm{CaT}}^2 \cdot  h_{\mathrm{CaT}} \cdot (v_{\mathrm{CT}} - \bar{E}_{\mathrm{Ca}}) \; - \; \bar{g}_{\mathrm{CaPQ}}\cdot m_{\mathrm{CaPQ}} \cdot (v_{\mathrm{CT}} - \bar{E}_{\mathrm{Ca}})}{2\cdot F \cdot \mathrm{\mathit{Vol}} \cdot  A} \right.\\
    &\left. \quad \quad \quad \quad \; \; - \; \frac{1}{k_v \cdot g_{\mathrm{max}} \cdot A} \frac{\mathrm{ca} - \bar{\mathrm{Ca}}_{eq}}{\tau_{store}} \right),
\end{aligned}
\end{equation}

\noindent where $\bar{g}_i=g_i / g_{\mathrm{max}}$ , $\bar{E}_i = E_i / k_v$ for $i \in \mathcal{I}_{\mathrm{CT}}$ (see Eq.~\eqref{EQN:CHOLT}), $\bar{I}_{\mathrm{App}} = I_{\mathrm{App}}/(k_v \cdot g_{\mathrm{max}})$, $\tau_{store}=t_{store}/k_{\tau}$, and $\bar{\mathrm{Ca}}_{eq} = \mathrm{Ca}_{eq} / k_{\mathrm{ca}}$ are dimensionless parameters. We define $R_p = k_{\tau} T_p$ for all gating variables $p$, $R_v := (g_{\mathrm{max}} \cdot k_{\tau})/c_m, R_{\mathrm{ca}}= A \cdot k_v \cdot g_{\mathrm{max}}$, and $A = \allowbreak 1/\max\{(2 \cdot F \cdot \mathit{Vol} \cdot k_{\mathrm{ca}})/k_{\tau},$ $ k_v \cdot g_{\mathrm{max}} \}$ = 0.2 $\mathrm{nA}^{-1}$. Then the functions on the right-hand side of the equations in System~\eqref{sys:appendix-ct-nondim1} remain $\mathcal{O}(1)$ over the relevant range of inputs. Values of dimensionless parameters $R_j$ for variables $j$ are summarized in~\Cref{TAB:ct-nondim}. We conclude that $x=[v_{\mathrm{CT}},m_{\mathrm{Na}}, m_{\mathrm{K}},m_{\mathrm{A}},h_{\mathrm{Na}},m_{\mathrm{CaPQ}}]^{\mathsf{T}}$ evolve on a fast time scale and form the fast subsystem, $y=[m_{\mathrm{CaT}}, h_{\mathrm{A}}]^{\mathsf{T}}$ the slow subsystem, and $z=[h_{\mathrm{CaT}},\mathrm{ca}]^{\mathsf{T}}$ the superslow subsystem. We choose the fast time scale as our reference time scale, and define
\[\varepsilon_1 := \max\{R_{m_\mathrm{CaT}}, R_{h_\mathrm{A}}\}, \; \varepsilon_2 := 1/\varepsilon_1 \cdot \max\{R_{h_\mathrm{CaT}}, R_{\mathrm{ca}}\}.\]

\noindent Substituting these expressions and consolidating our notation for the right-hand side functions in the model, we obtain a system of the form~\eqref{sys:ct-nondim} by rewriting the dimensionless system 

\begin{figure}[H]
\centering
\includegraphics{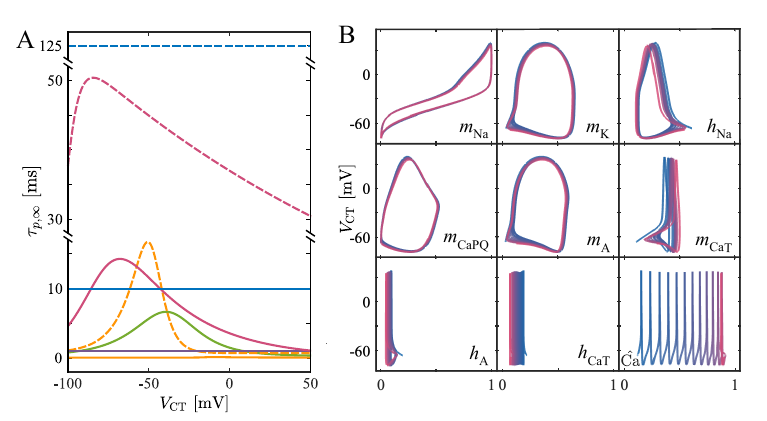}
\caption{\label{fig:cholt-pif_timescales} \small Timescales and trajectories during oscillations in the CT model response to the PIF protocol. 
\textbf{A}: Time constants for activation $m_{\mathrm{Na}}$ (solid) and inactivation $h_{\mathrm{Na}}$ (dashed) gating variables depend nonlinearly on voltage: $\mathrm{Na}$ (orange), $\mathrm{K}$ (green), P/Q-type $\mathrm{Ca}$ (purple), A-type $\mathrm{K}$ (pink), T-type $\mathrm{Ca}$ (blue). 
\textbf{B}: Trajectories of CT model variables, indicated in each panel, during PIF response, visualized against dimensionalized $V_{\mathrm{CT}}$, colored from $t=290$ (blue) to $t=600$ (pink). The subpanels for $m_{\mathrm{Na}}, m_{\mathrm{K}}, h_{\mathrm{Na}}, m_{\mathrm{CaPQ}},$ and $m_{\mathrm{A}}$ show periodic orbits with minimal drift. Those for $m_{\mathrm{CaT}}$ and $h_A$ also show oscillations but with relatively more drift. The panels for $h_{\mathrm{CaT}}, \mathrm{\hat{Ca}}$, in contrast, show a systematic drift (to obtain $\mathrm{\hat{Ca}}$ for visualization, we transformed values of $\mathrm{Ca}$ (nM) using the linear function $\hat{\mathrm{Ca}} = (\mathrm{Ca}-100) / 30$).
}
\end{figure}

\begin{equation}\label{sys:appendix-ct-nondim2}
\begin{aligned}
    \frac{dv}{d\tau} &= \; R_v \Tilde{f}_1(x,y,z) = f_1(x,y,z),\\
    \frac{dm_{\mathrm{Na}}}{d\tau} &= \; R_{m_{\mathrm{Na}}} \Tilde{f}_2(x,y,z) = f_2(x,y,z),\\
    \frac{dm_{\mathrm{K}}}{d\tau} &= \; R_{m_{\mathrm{K}}} \Tilde{f}_3(x,y,z) = f_3(x,y,z),\\
    \frac{dm_{\mathrm{A}}}{d\tau} &= \; R_{m_{\mathrm{A}}} \Tilde{f}_4(x,y,z) = f_4(x,y,z),\\
    \frac{dh_{\mathrm{Na}}}{d\tau} &= \; R_{h_{\mathrm{Na}}} \Tilde{f}_5(x,y,z) = f_5(x,y,z),\\
    \frac{dm_{\mathrm{CaPQ}}}{d\tau} &= \; R_{m_{\mathrm{CaPQ}}} \Tilde{f}_6(x,y,z) = f_6(x,y,z),\\
    \frac{dm_{\mathrm{CaT}}}{d\tau} &= \; R_{m_{\mathrm{CaT}}} \Tilde{g}_1(x,y) = \varepsilon_1 g_1(x,y),\\
    \frac{dh_{\mathrm{A}}}{d\tau} &= \; R_{h_{\mathrm{A}}} \Tilde{g}_2(x,y) = \varepsilon_1 g_2(x,y),\\
    \frac{dh_{\mathrm{CaT}}}{d\tau} &= \; R_{h_{\mathrm{CaT}}} \Tilde{h}_1(x,y,z) = \varepsilon_1 \varepsilon_2 h_1(x,y,z),\\
    \frac{d\mathrm{ca}}{d\tau} &= \; R_{\mathrm{ca}} \Tilde{h}_2(x,y,z) = \varepsilon_1 \varepsilon_2 h_2(x,y,z).
\end{aligned}
\end{equation}

In the PIF response by the CT model in~\Cref{sec:ct-pif}, the analysis leverages the same three-timescale decomposition with fast $x$, slow $y$ and superslow $z$ subsystems (see~\Cref{TAB:ct-nondim}) based on the same time constants for the analysis of the CT model response in the PIR (\Cref{sec:ct-pir}) protocol.

\subsection{Two timescales within the C PPN model}\label{sec:appendix-c-nondim}

For the C model during the post-inhibitory delay response, a similar nondimensionalization procedure to the previous section~(\Cref{sec:appendix:ct-nondim}) reveals the time constants listed in~\Cref{tab:c-timeconstants}. We identify two timescales within the model, casting it as a fast-slow system of the form~\eqref{sys:nc-ramp-fasttime} with $x = [v_{\mathrm{C}}, m_{\mathrm{Na}}, m_{\mathrm{K}}, \allowbreak m_{\mathrm{CaPQ}}, h_{\mathrm{Na}}, m_{\mathrm{A}}]^{\mathsf{T}}$ and $y = [h_{\mathrm{A}},\mathrm{ca}]^{\mathsf{T}}$ as fast and slow subsystems, respectively. Functions $f,g$ and parameter $\varepsilon_1$ are defined analogously as in System~\eqref{sys:appendix-ct-nondim2} in~\Cref{sec:appendix:ct-nondim}.

\begin{table}[ht!]\centering
\caption{\label{tab:c-timeconstants} Time constants $R_x$ and timescale classifications for nondimensionalization of C model in response to the depolarizing delay protocol in~\Cref{sec:c-delay}.}
\renewcommand{\arraystretch}{1.3}
$\begin{array}{l|llllllll}
x & v_{\mathrm{C}} & m_{\mathrm{Na}} & m_{\mathrm{K}} & m_{\mathrm{CaPQ}} & h_{\mathrm{Na}} & m_{\mathrm{A}} & h_{\mathrm{A}} & \mathrm{ca}
\\
\hline
R_x \; \left[10^{-3}\right] & 50 & 20 & 1.01 & 1 & 0.72 & 0.48 & 0.07 & 0.001 \\
\mathrm{timescale} & \mathrm{fast} & \mathrm{fast} & \mathrm{fast} & \mathrm{fast} & \mathrm{slow} & \mathrm{slow} & \mathrm{slow} & \mathrm{slow}
\end{array}$
\end{table}

\subsection{Multiple timescales within the NC PPN model}\label{sec:appendix-nc-nondim}

A similar procedure to that in~\Cref{sec:appendix:ct-nondim}, applied to the NC model during the depolarizing ramp, SDP, and PIF protocols, reveals the nondimensionalization parameters listed in~\Cref{tab:nc-timeconstants}. For each applied current protocol, we identify two timescales within the model cell and analyze the model as a fast-slow system of the form~\eqref{sys:nc-ramp-fasttime} where $f,g$ and $\varepsilon_1$ are defined analogously as in~\eqref{sys:appendix-ct-nondim2} in~\Cref{sec:appendix:ct-nondim}. During the response to the ramping protocol where $I_{\mathrm{Na}}$ channels are blocked, we identify fast and slow subsystems $x = [v_{\mathrm{NC}}, m_{\mathrm{Na}},m_{\mathrm{K}},h_{\mathrm{Na}},m_{\mathrm{CaPQ}}, m_{\mathrm{CaT}}]^{\mathsf{T}}$ and $y = [h_{\mathrm{CaT}},\mathrm{ca}]^{\mathsf{T}}$, respectively. For the SDP protocol, $x = [v_{\mathrm{NC}}, m_{\mathrm{Na}},m_{\mathrm{K}},h_{\mathrm{Na}},m_{\mathrm{CaPQ}}]^{\mathsf{T}}$ form the fast variables and $y = [m_{\mathrm{CaT}},h_{\mathrm{CaT}},\mathrm{ca}]^{\mathsf{T}}$ represent the slow variables. 

For the model activity under the PIF protocol, our analysis requires two distinct decompositions into fast-slow subsystems. During the slow depolarization of the cell before spiking, we use a narrow range of $v_{\mathrm{NC}}$ corresponding to the slow depolarization phase to calculate the second set of time constants in~\Cref{tab:nc-timeconstants}, where the fast variables are $\hat{x}=[v_{\mathrm{NC}}, m_{\mathrm{Na}},m_{\mathrm{K}},h_{\mathrm{Na}}]^{\mathsf{T}}$ and slow variables are $\hat{y}=[m_{\mathrm{CaPQ}}, m_{\mathrm{CaT}},h_{\mathrm{CaT}},\mathrm{ca}]^{\mathsf{T}}$. The remainder of the activity is analyzed using another decomposition, where $x=[v_{\mathrm{NC}}, m_{\mathrm{Na}},m_{\mathrm{CaPQ}}]^{\mathsf{T}}$ are fast and $y=[m_{\mathrm{K}},h_{\mathrm{Na}}, m_{\mathrm{CaT}},h_{\mathrm{CaT}},\allowbreak \mathrm{ca}]^{\mathsf{T}}$ are slow. To accommodate this variability, when we present the timescale decomposition in~\Cref{sec:nc-ct-pif}, we split off $m_{\mathrm{K}}$ and $h_{\mathrm{Na}}$ into their own slow class and label the remaining variables as fast or superslow.

\begin{table}[ht!]\centering
\caption{\label{tab:nc-timeconstants} Time constants $R_x$ and timescale classifications for nondimensionalization of NC model in response to the depolarizing ramp, SDP and PIF protocols. See~\Cref{sec:nc-pif} for details on timescale classifications for PIF.}
\renewcommand{\arraystretch}{1.3}

\begin{tabular}{l||l|llllllll}
Protocol & $x$ & $v_{\mathrm{NC}}$ & $m_{\mathrm{Na}}$ & $m_{\mathrm{K}}$ & $h_{\mathrm{Na}}$ & $m_{\mathrm{CaPQ}}$ & $m_{\mathrm{CaT}}$ & $h_{\mathrm{CaT}}$ & $\mathrm{ca}$ \\
\hline
\multirow{2}{*}{$\mathrm{Ramp.} \; (\ref{sec:nc-ramp})$}
& $R_x \left[10^{-3}\right]$ & 50 & $-$ & 0.29 & $-$ & 1 & 0.1 & 0.008 & 0.001 \\
& $\mathrm{timescale}$ & fast & $-$ & fast & $-$ & fast & fast & slow & slow \\
\hline
\multirow{2}{*}{$\mathrm{SDP} \; (\ref{sec:nc-sdp})$}
& $R_x \left[10^{-3}\right]$ & 50 & 20 & 1.78 & 1.43 & 1 & 0.1 & 0.008 & 0.001 \\
& $\mathrm{timescale}$ & fast & fast & fast & fast & fast & slow & slow & slow \\
\hline

\multirow{4}{*}{$\mathrm{PIF} \; (\ref{sec:nc-ct-pif})$}
& $R_x \left[10^{-3}\right]$ & 50 & 20 & 0.32 & 0.14 & 1 & 0.1 & 0.008 & 0.001 \\
& timescale-1 & & & & & & & & \\
\cline{2-10}
& $R_x \left[10^{-3}\right]$ & 50 & 20 & 1.26 & 1.43 & 1 & 0.1 & 0.008 & 0.001 \\
& timescale-2 & & & & & & & \\
\hline
\end{tabular}
\end{table}

\section*{Acknowledgments}
\noindent This work was partially supported by NIH awards R01NS125814 and R01DA059993, both part of the CRCNS program. The authors thank Rebekah Evans, Aryn Gittis, and their lab members for extremely helpful discussions about PPN neurons. A.T. thanks Matteo Martin for insightful discussions on multiple timescale analysis and guidance with XPPAUT.

\bibliographystyle{abbrv}
\bibliography{references}

\end{document}